\documentclass{article}

\usepackage{varioref}
\usepackage{srcltx}
\input{xy}
\xyoption{all}

\usepackage{verbatim}
\usepackage{varioref}
\usepackage{amsfonts,
enumerate,
amsthm,
amsmath,latexsym,amscd,amssymb,doc}
\usepackage[latin1]{inputenc}
\usepackage{graphicx}
\usepackage{hyperref}

\usepackage{cyrillic}
\newcommand{\cyrrm}{\fontencoding{OT2}\selectfont\textcyrup}

\newcommand{\Beilinson}{{\cyrrm{B}}} 

\usepackage{xr}
		\newcommand{\category}[1]{\mathbf{#1}}
	
		\newcommand{\Com}{\category{Com}} 
    \newcommand{\D}{\category {D}} 

    \newcommand{\Ab}{\category {Ab}} 
    
    \newcommand{\Sm}{\category{Sm}}

		\newcommand{\Spt}{\category{Spt}} 
		\newcommand{\PSh}{\category{PSh}} 

    \newcommand{\DM}{\category {DM}}
    \newcommand{\DMBGL}{\DM_{\BGL}}
    \newcommand{\DMBei}{\DM_{\Beilinson}}

\newcommand{\Mod}{\category {Mod}} 

\newcommand{\SH}{\category{SH}} 
\newcommand{\HoC}{\category{Ho}} 
\newcommand{\HoCsect}{\HoC_{\mathrm{sect}, \bullet}} 

	    \newtheoremstyle{Normal}
  {}
  {}
  {}
  {}
  {\bfseries}
  {.}
  { }
  {}

    \theoremstyle{Normal} 

    \newtheorem{Defi}{Definition}[section]

    \newtheorem{Conj}[Defi]{Conjecture}
    \newtheorem{Bem}[Defi]{Remark}
    \newtheorem{Bsp}[Defi]{Example}
     
    \newtheorem{Axio}[Defi]{Axiom}
    \newtheorem{Ques}[Defi]{Question}
    \newtheorem{Comp}[Defi]{Complements}
    
		\theoremstyle{remark}
	
    \theoremstyle{plain} 

    \newtheorem{Satz}[Defi]{Proposition}
    \newtheorem{DefiTheo}[Defi]{Definition and Theorem}
    \newtheorem{Theo}[Defi]{Theorem}
    \newtheorem{Folg}[Defi]{Corollary}
    \newtheorem{Lemm}[Defi]{Lemma}
    \newtheorem{DefiLemm}[Defi]{Definition and Lemma}

\newcommand{\refit}[1]{(\ref{item_#1})}

\newcommand{\refsect}[1]{Section \ref{sect_#1}}
\newcommand{\conj}{\begin{Conj}} 			\newcommand{\xconj}{\end{Conj}}											    
\newcommand{\ques}{\begin{Ques}} 			\newcommand{\xques}{\end{Ques}}											    
\newcommand{\axio}{\begin{Axio}} 			\newcommand{\xaxio}{\end{Axio}}											    
\newcommand{\bem}{\begin{Bem}} 			\newcommand{\xbem}{\end{Bem}}											    
\newcommand{\rema}{\begin{Bem}} 			\newcommand{\xrema}{\end{Bem}}											    \newcommand{\refre}[1]{Remark \ref{rema_#1}}
\newcommand{\defi}{\begin{Defi}} 			\newcommand{\xdefi}{\end{Defi}}										\newcommand{\refde}[1]{Definition \ref{defi_#1}}
\newcommand{\defitheo}{\begin{DefiTheo}} \newcommand{\xdefitheo}{\end{DefiTheo}} 
\newcommand{\defilemm}{\begin{DefiLemm}} \newcommand{\xdefilemm}{\end{DefiLemm}} \newcommand{\refdele}[1]{Definition and Lemma \ref{defilemm_#1}}

\newcommand{\lemm}{\begin{Lemm}}			\newcommand{\xlemm}{\end{Lemm}}											\newcommand{\refle}[1]{Lemma \ref{lemm_#1}}
\newcommand{\comp}{\begin{Comp}}			\newcommand{\xcomp}{\end{Comp}}											

\newcommand{\satz}{\begin{Satz}}			\newcommand{\xsatz}{\end{Satz}}										

\newcommand{\prop}{\begin{Satz}}			\newcommand{\xprop}{\end{Satz}}										\newcommand{\refpr}[1]{Proposition \ref{prop_#1}}
\newcommand{\theo}{\begin{Theo}}			\newcommand{\xtheo}{\end{Theo}}											\newcommand{\refth}[1]{Theorem \ref{theo_#1}}
\newcommand{\bsp}{\begin{Bsp}}				\newcommand{\xbsp}{\end{Bsp}}												
\newcommand{\exam}{\begin{Bsp}}				\newcommand{\xexam}{\end{Bsp}}												\newcommand{\refex}[1]{Example \ref{exam_#1}}
\newcommand{\folg}{\begin{Folg}}				\newcommand{\xfolg}{\end{Folg}}
				 			
\newcommand{\cor}{\begin{Folg}}				\newcommand{\xcor}{\end{Folg}}
\newcommand{\mycomment}{}

\newcommand{\eqnarra}{\begin{eqnarray}}				\newcommand{\xeqnarra}{\end{eqnarray}}
\newcommand{\eqnarr}{\begin{eqnarray*}}				\newcommand{\xeqnarr}{\end{eqnarray*}}

\newcommand{\eqn}{\begin{equation}} 		\newcommand{\xeqn}{\end{equation}}
\newcommand{\refeq}[1]{(\ref{eqn_#1})}

  \newcommand{\gl}{\index}

  \newcommand{\Dph}[1]{\emph{#1}}
  \newcommand{\Def}[1]{\Dph{#1}\gl{#1}}

\newcommand{\mylabel}[1]{\label{#1}}

\newcommand{\status}[1]{}
\newcommand{\explain}[1]{\stackrel{\text{\refeq{#1}}}{=}}

\newcommand{\BB}{\mathrm{B}}
\newcommand{\Ho}{\category{Ho}}
\newcommand{\OF}{{\mathcal{O}_F}}

\newcommand{\twi}[1]{\{#1\}}

\renewcommand{\log}{\mathrm{log}\,}
\newcommand{\C}{\mathcal{C}}

\newcommand{\R}{\mathrm{R}} 
\newcommand{\RG}{{\R} {\Gamma}}
\newcommand{\HD}{\H_\mathrm{D}} 
\newcommand{\HDel}[1]{\H_{\mathrm{D}, #1}} 

\newcommand{\lr}{{\longrightarrow}}
\renewcommand{\r}{\rightarrow}

\renewcommand{\t}{{\otimes}}

\newcommand{\del}{\partial}

\newcommand{\A}[1][1]{\mathbb{A}^{#1}}
\renewcommand{\P}[1][1]{\mathbb P^{#1}}

\newcommand{\Q}{\mathbb{Q}}

\newcommand{\CC}{\mathbb{C}}
\newcommand{\RR}{\mathbb{R}}
\newcommand{\Z}{\mathbb{Z}}

\renewcommand{\H}{\mathrm{H}}

\newcommand{\x}{{\times}}
\newcommand{\ol}[1]{{\overline{#1}}}

\newcommand{\Beweis}{{\normalfont} \textbf{Proof}}

\newcommand{\lc}{\textit{loc.~cit.}}

\newcommand{\lcs}{\textit{loc.~cit.}\ }

\newcommand{\Spec}{\mathrm{Spec} \; \! }
\newcommand{\SpecOF}{{\Spec} {\OF}}

\newcommand{\SpecZ}{{\Spec} {\Z}}
\newcommand{\SpecQ}{{\Spec} {\Q}}

\newcommand{\Fr}{\operatorname{Fr}} 
\newcommand{\Td}{\operatorname{Td}} 

\newcommand{\op}{\mathrm{op}}

\newcommand{\red}{\mathrm{red}}

\newcommand{\GmS}{\mathbb{G}_{m,S}}

\newcommand{\gr}{\operatorname{gr}} 

\newcommand{\tr}{\mathrm{tr}}
\newcommand{\ch}{\operatorname{ch}} 
\newcommand{\dR}{\mathrm{dR}}

\newcommand{\Hom}{\mathrm{Hom}}
\newcommand{\End}{\mathrm{End}}
\newcommand{\IHom}{\underline{\Hom}}
\newcommand{\id}{\mathrm{id}}
\renewcommand{\gr}{\operatorname{gr}} 

\newcommand{\im}{\operatorname{im}} 

\newcommand{\M}{{\operatorname{M}}} 

\newcommand{\coker}{\operatorname{coker}}

\newcommand{\cone}{\operatorname{cone}}
\newcommand{\Tot}{\mathrm{Tot}}

\newcommand{\loc}[2]{[}

\newcommand{\pr}{\begin{proof}[\Beweis: ]}
\newcommand{\pf}{\pr}
\newcommand{\xpf}{\end{proof}}

\newcommand{\SiginfP}{\Sigma^\infty_{\P}} 

\theoremstyle{definition}

\theoremstyle{definition}

\theoremstyle{definition}

\theoremstyle{definition}

\theoremstyle{definition}

\theoremstyle{definition}

\theoremstyle{definition}

\theoremstyle{definition}

\theoremstyle{definition}

\theoremstyle{definition}

\theoremstyle{definition}

\theoremstyle{definition}

\theoremstyle{definition}

\theoremstyle{definition}

\newcommand{\HB}{\mathrm{H}_{\Beilinson}}
\newcommand{\HBei}[1]{\mathrm{H}_{\Beilinson, #1}}
\newcommand{\HBhat}{\widehat{\HB}}
\newcommand{\HBeihat}[1]{\widehat{\HBei{#1}}}
\newcommand{\HBhatS}{\HBeihat{S}}

\newcommand{\HBS}{\mathrm{H}_{{\Beilinson}, S }}

\newcommand{\Hhat}{\widehat{\mathrm{H}}}
\newcommand{\Khat}{\widehat{\mathrm{K}}}

\newcommand{\hofib}{\operatorname{hofib}} 

\newcommand{\Pic}{\mathrm{Pic}}

\newcommand{\cK}{\mathcal{K}}
\newcommand{\cO}{\mathcal{O}}

\newcommand{\cDz}{\H_\Deligne} 

\newcommand{\Deligne}{\mathrm{D}}

\newcommand{\BGL}{\mathrm{BGL}}
\newcommand{\BGLhat}{\widehat{\BGL}}
\newcommand{\BGLhatS}{\widehat{\BGL_S}}

\numberwithin{equation}{section}

\makeatother

\begin{document}

\author{Andreas Holmstrom\footnote{I.H.E.S. Le Bois-Marie, 35 Route de Chartres, F-91440 Bures-sur-Yvette, France, \href{mailto:andreas.holmstrom@gmail.com}{andreas.holmstrom@gmail.com}} \ and Jakob Scholbach\footnote{Universit{\"a}t M{\"u}nster, Mathematisches Institut, Einsteinstr. 62, D-48149 M{\"u}nster, Germany}
}

\title{Arakelov motivic cohomology I}


\maketitle


\begin{abstract}
This paper introduces a new cohomology theory for schemes of finite type over an arithmetic ring. The main motivation for this Arakelov-theoretic version of motivic cohomology is the conjecture on special values of $L$-functions and zeta functions formulated by the second author \cite{Scholbach:specialL}. Taking advantage of the six functors formalism in motivic stable homotopy theory, we establish a number of formal properties, including pullbacks for arbitrary morphisms, pushforwards for projective morphisms between regular schemes, localization sequences, $h$-descent. We round off the picture with a purity result and a higher arithmetic Riemann-Roch theorem.

In the second part of the paper, we relate Arakelov motivic cohomology to classical constructions such as arithmetic $K$ and Chow groups and the height pairing.
\end{abstract}

\section{Introduction}
For varieties over finite fields, we have very good cohomological tools for understanding the associated zeta functions. These tools include $\ell$-adic cohomology, explaining the functional equation and the Riemann hypothesis, and Weil-\'etale cohomology, which allows for precise conjectures and some partial results regarding the ``special values'', i.e. the vanishing orders and leading Taylor coefficients at integer values. The conjectural picture for zeta functions of schemes $X$ of finite type over $\SpecZ$ is less complete. Deninger envisioned a cohomology theory explaining the Riemann hypothesis, and Flach and Morin have developed the Weil-\'etale cohomology  describing special values of zeta functions of regular projective schemes over $\Z$ at $s=0$ \cite{Deninger:Motivic, FlachMorin:Weil, Morin:Zeta}.

In \cite{Scholbach:specialL}, the second author proposed a new conjecture, which describes the special values of all zeta functions and $L$-functions of geometric origin, up to a rational factor. It is essentially a unification of classical conjectures of Beilinson, Soul\'e and Tate, formulated in terms of  the recent Cisinski-D\'eglise theory of triangulated categories of motives over $\Z$. This conjecture is formulated in terms of a new cohomology theory for schemes of finite type over $\Z$. The purpose of this paper is to construct this cohomology theory and establish many of its properties.

This cohomology theory, we call it \emph{Arakelov motivic cohomology}, is related to motivic cohomology, roughly in the same way as arithmetic Chow groups relate to ordinary Chow groups or as arithmetic $K$-theory relates to algebraic $K$-theory. The key principle for cohomology theories of this type has always been to connect some algebraic data, such as the algebraic $K$-theory, with an analytical piece of information, chiefly Deligne cohomology, in the sense of long exact sequences featuring the Beilinson regulator map between the two and a third kind of groups measuring the failure of the regulator to be an isomorphism. This was suggested by Deligne and Soul\'e in the 80s. Beilinson also expressed the idea that the "boundary" of an algebraic cycle on a scheme over $\Z$ should be a Deligne cohomology class \cite{Beilinson:Height}. Gillet, Roessler, and Soul\'e then started developing a theory of arithmetic Chow groups \cite{GilletSoule:CharacteristicI, GilletSoule:CharacteristicII, GilletSoule:Arithmetic, Soule:Arakelov}, arithmetic $K_0$-theory and an arithmetic Riemann-Roch theorem \cite{Roessler:Adams,GilletRoesslerSoule}. Burgos and Wang \cite{Burgos:Dolbeault, Burgos:Arithmetic, BurgosWang} extended some of this to not necessarily projective schemes and gave an explicit representation of the Beilinson regulator. More recently, Goncharov gave a candidate for higher arithmetic Chow groups for complex varieties, Takeda developed higher arithmetic $K$-theory, while Burgos and Feliu constructed higher arithmetic Chow groups for varieties over arithmetic fields \cite{Goncharov:Regulators, Takeda:Higher, BurgosFeliu:Higher}. The analogous amalgamation of topological $K$-theory and Deligne cohomology of smooth manifolds is known as smooth $K$-theory \cite{BunkeSchick:Smooth}.

In a nutshell, these constructions proceed by representing the regulator as a map of appropriate complexes. Then one defines, say, arithmetic $K$-theory to be the cohomology of the cone of this map. Doing so, however, requires a good command of the necessary complexes, which so far prevented extending higher arithmetic Chow groups to schemes over $\Z$ and also requires to manually construct homotopies whenever a geometric construction is to be done, for example the pushforward. The idea of this work is to both overcome these hurdles and enhance the scope of these techniques by introducing a \emph{spectrum}, i.e., an object in the stable homotopy category of schemes, representing the sought cohomology theory.

This paper can be summarized as follows:
let $S$ be a regular scheme of finite type over a number field $F$, a number ring $\OF$, $\RR$, or $\CC$. In the stable homotopy category $\SH(S)$ (cf.\ \refsect{SH}) there is a ring spectrum $\cDz$ representing Deligne cohomology with real coefficients of smooth schemes $X / S$ (\refth{Delignespectrum}). 
We \emph{define} (cf.\ \ref{defi_widehatA}) the Arakelov motivic cohomology spectrum $\HBhat$ as the homotopy fiber of the map
$$\HB \stackrel{\id \wedge 1_{\HD}}\lr \HB \wedge \HD.$$
Here, $\HB$ is Riou's spectrum representing the Adams eigenspaces in algebraic $K$-theory (tensored by $\Q$). \'Etale descent for $\HD$ implies that the canonical map $\HD \r \HB \wedge \HD$ is an isomorphism (\ref{theo_Delignespectrum}), so there is a distinguished triangle
$$\HBhat \r \HB \r \HD \r \HBhat[1].$$
We define \Def{Arakelov motivic cohomology} to be the theory represented by this spectrum, that is to say
$$\Hhat^n(M, p) := \Hom_{\SH(S)_\Q}(M, \HBhat(p)[n])$$
for any $M \in \SH(S)$. Thus, there is a long exact sequence involving Arakelov motivic cohomology, motivic cohomology and Deligne cohomology (\refth{longexactArakelov}). Moreover, Arakelov motivic cohomology shares the structural properties known for motivic cohomology, for example a projective bundle formula, a localization sequence, and $h$-descent (\refth{immediate}). It also has the expected \emph{functoriality}: pullback for arbitrary morphisms of schemes (or motives, \refle{pullback}) and pushforward along projective maps between regular schemes (\refdele{pushforward}). All of this can be modified by replacing $\HB$ by $\BGL$, the spectrum representing algebraic $K$-theory. The resulting Arakelov version is denoted $\BGLhat$ and the cohomology theory so obtained is denoted $\Hhat^n(M)$.

We extend the motivic Riemann-Roch theorem given by Riou to arbitrary projective maps between regular schemes (\refth{RiemannRoch})---a statement that is of independent interest. We deduce a \emph{higher arithmetic Riemann-Roch theorem} (\refth{HARR}) for the cohomology theories $\Hhat^*(M, -)$ vs.\ $\Hhat^*(M)$. It applies to smooth projective morphisms and for projective morphisms between schemes that are smooth over the base.\\

In the second part of this paper, we will show how to relate the homotopy-theoretic construction of Arakelov motivic cohomology to the classical definitions of arithmetic $K$- and Chow groups. For example, the arithmetic $K_0$-groups $\Khat_0^T(X)$ defined by Gillet and Soul\'e \cite[Section 6]{GilletSoule:CharacteristicII} for a regular projective variety $X$ (over a base $S$ as above) sit in an exact sequence
$$K_1(X) \r \oplus_{p \geq 0} A^{p,p}(X) / (\im \partial + \im \overline \partial) \r \Khat^T_0(X) \r K_0(X) \r 0,$$
where $A^{p,p}(X)$ is the group of real-valued $(p,p)$-forms $\omega$ on $X(\CC)$ such that $\Fr_\infty^* \omega = (-1)^p \omega$. The full arithmetic $K$-groups $\Khat^T_0(X)$ are not homotopy invariant and can therefore not be adressed using $\A$-homotopy theory. Instead, we consider the subgroup
$$\Khat_0(X) := \ker \left (\ch: \Khat^T_0(X) \r \oplus_{p \geq 0} A^{p,p}(X) \right ).$$
For smooth schemes $X / S$, we show a canonical isomorphism
\eqn \mylabel{eqn_KhatintroI}
\Hhat^0(\M(X)) \cong \Khat_0(X)
\xeqn
and similarly for higher arithmetic $K$-theory, as defined by Takeda. The homotopy-theoretic approach taken yields a considerable simplification since it is no longer necessary to construct explicit homotopies between the complexes representing arithmetic $K$-groups, say. For example, the Adams operations on $\Khat_i(X)$ defined by Feliu \cite{Feliu:Adams} were not known to induce a decomposition $\Khat_*(X)_\Q \cong \oplus_p \Khat_*(X)^{(p)}_\Q$. Using that the isomorphism \refeq{KhatintroI} is compatible with Adams operations, this statement follows from the essentially formal analogue for $\Hhat^*$. Moreover, \refeq{KhatintroI} is shown to be compatible with the pushforwards on both sides in an important case. This implies that the height pairing on a smooth projective scheme $X / S$, $S \subset \Spec \Z$, is expressible in terms of the natural pairing of motivic homology and Arakelov motivic cohomology of the motive of $X$. According to the second author's conjecture, the $L$-values of schemes (or motives) over $\Z$ are given by the determinant of this pairing.
%

\hspace{0cm}

It is a pleasure to thank Denis-Charles Cisinski and Fr\'ed\'eric D\'eglise for a number of enlightening conversations. We also thank the referee for a number of helpful comments. The first-named author also wishes to thank Tony Scholl and Peter Arndt. The second-named author gratefully acknowledges the hospitality of Universit\'e Paris 13, where part of this work was done.

%
%

\section{Preliminaries} \mylabel{sect_preliminaries}
In this section, we provide the motivic framework that we are going to work with in Sections \ref{sect_Deligne} and \ref{sect_Arakelov}: we recall the construction of the stable homotopy category $\SH(S)$ and some properties of the Cisinski-D\'eglise triangulated category of motives. In \refsect{RiemannRoch}, we generalize Riou's formulation of the Riemann-Roch theorem to regular projective morphisms. This will then be used to derive a higher arithmetic Riemann-Roch theorem (\ref{theo_HARR}). Finally, we recall the definition and basic properties of Deligne cohomology that are needed in \refsect{Deligne} to construct a spectrum representing Deligne cohomology.

\subsection{The stable homotopy category} \mylabel{sect_SH}

This section sets the notation and recalls some results pertaining to the homotopy theory of schemes due to Morel and Voevodsky \cite{MorelVoevodsky:A1}.

Let $S$ be a Noetherian scheme. We only use schemes which are of finite type over $\Z$, $\Q$, or $\RR$. Unless explicitly mentioned otherwise, all morphisms of schemes are understood to be separated and of finite type. Let $\Sm / S$ be the category of smooth schemes over $S$. The category of presheaves of pointed sets on this category is denoted $\PSh_\bullet := \PSh_\bullet(\Sm / S)$. We often regard a scheme $X \in \Sm /S$ as the presheaf (of sets) represented by $X$, and we write $X_+ := X \sqcup \{ * \}$ for the associated pointed version. The projective line $\P_S$ is always viewed as pointed by $\infty$. The prefix $\Delta^\op -$ indicates simplicial objects in a category. The simplicial $n$-sphere is denoted $S^n$, this should not cause confusion with the base scheme $S$.

We consider the \emph{pointwise} and the \emph{motivic model structure} on the category $\Delta^\op (\PSh_\bullet)$ \cite[Section 1.1.]{Jardine:Motivic}. The latter is obtained by considering objects that are local with respect to projections $U \x \A \r U$  and the Nisnevich topology. The corresponding homotopy categories will be denoted by $\HoCsect$ and $\Ho_\bullet$, respectively. The identity functor is a Quillen adjunction with respect to these two model structures.

The category $\Spt := \Spt^{\P}(\Delta^\op \PSh_\bullet(\Sm/S))$ consists of symmetric $\P_S$-spectra, that is, sequences $E = (E_n)_{n \geq 0}$ of simplicial presheaves which are equipped with an action of the symmetric group $S_n$ and bonding maps $\P \wedge E_n \r E_{n+1}$ such that $(\P)^{\wedge m} \wedge E_n \r E_{n+m}$ is $S_n \x S_m$-equivariant (and the obvious morphisms). The functor $\SiginfP : \Delta^\op(\PSh_\bullet) \ni F \mapsto \left ((\P)^{\wedge n} \wedge F \right)_{n \geq 0}$ (bonding maps are identity maps, $S_n$ acts by permuting the factors $\P$) is left adjoint to $\Omega^\infty : (E_n) \mapsto E_0$. Often, we will not distinguish between a simplicial presheaf $F$ and $\SiginfP(F)$.

The category $\Spt$ is endowed with the \Def{stable model structure} \cite[Theorems 2.9, 4.15]{Jardine:Motivic}. The corresponding homotopy category is denoted $\SH$ (or $\SH(S)$) and referred to as the \emph{stable homotopy category} of smooth schemes over $S$. The pair $(\SiginfP, \Omega^\infty)$ is a  Quillen adjunction with respect to the motivic model structure on $\Delta^\op \PSh_\bullet$ and the stable model structures on $\Spt$. We sum up this discussion by saying that there are adjunctions of homotopy categories
\eqn \mylabel{eqn_SH}
\HoCsect \leftrightarrows \Ho_{\bullet} \leftrightarrows \SH.
\xeqn

The stable homotopy categories are triangulated categories. We will use both the notation $M[p]$ and $M \wedge (S^1)^{\wedge p}$, $p \in \Z$ for the shift functor. Moreover, \mylabel{smashinginvertible} in $\HoC(S)$, there is an isomorphism $\P_S \cong S^1 \wedge (\GmS, 1)$. Thus, in $\SH(S)$, wedging with $\GmS$ is invertible, as well, and we write $M(p)$ for $M \wedge (\GmS)^{\wedge p}[-p]$, $p \in \Z$ for the \emph{Tate twist}. For brevity, we also put
$$M\twi{p} := M(p)[2p].$$

For any triangulated, compactly generated category $\C$ that is closed under coproducts, we let $\C_\Q$ be the full triangulated subcategory of those objects $Y$ such that $\Hom_{\C}(-, Y)$ is a $\Q$-vector space. The inclusion $i: \C_\Q \subset \C$ has a right adjoint which will be denoted by $(-)_\Q$. The natural map $\Hom_{\C}(X, Y) \t \Q \r \Hom_{\C}(X, i(Y_\Q)) = \Hom_{\C_\Q}(X_\Q, Y_\Q)$ is an isomorphism if $X$ is compact, see e.g. \cite[Appendix A.2]{Riou:Operations}. In particular, we will use $\SH(S)_\Q$. \mylabel{rational_coeff} Wherever convenient, we use the equivalence of this category with $\D_{\A}(S, \Q)$, the homotopy category of symmetric $\P$-spectra of complexes of Nisnevich sheaves of $\Q$-vector spaces (with the Tate twist and $\A$ inverted) \cite[5.3.22, 5.3.37]{CisinskiDeglise:Triangulated}.

Given a morphism $f: T \r S$, the stable homotopy categories are connected by adjunctions:
\eqn
\mylabel{eqn_adjstar}
f^* : \SH(S) \rightleftarrows \SH (T) : f_*,
\xeqn
\eqn
\mylabel{eqn_adjshriek}
f_! : \SH(T) \rightleftarrows \SH (S) : f^!,
\xeqn
\eqn
\mylabel{eqn_adjsharp}
f_\sharp : \SH(T) \rightleftarrows \SH (S) : f^*.
\xeqn
For the last adjunction, $f$ is required to be smooth. \refeq{adjstar} also applies to morphisms which are not necessarily of finite type (\cite[Scholie 1.4.2]{Ayoub:Six1}, see also \cite[1.1.11, 1.1.13; 2.4.4., 2.4.10]{CisinskiDeglise:Triangulated}).

\subsection{Beilinson motives} \mylabel{sect_motives}

Let $S$ be a Noetherian scheme of finite dimension. The key to Beilinson motives (in the sense of Cisinski and D\'eglise) is the motivic cohomology spectrum $\HBS$ due to Riou \cite[IV.46, IV.72]{Riou:Operations}. 
There is an object $\BGL_S \in \SH(S)$ representing algebraic $K$-theory in the sense that
\eqn \mylabel{eqn_BGL_vs_SH}
\Hom_{\SH(S)}(S^n \wedge \SiginfP X_+, \BGL_S) = K_n(X)
\xeqn
for any regular scheme $S$ and any smooth scheme $X / S$, functorially (with respect to pullback) in $X$. The $\Q$-localization $\BGL_{S, \Q}$ decomposes as
$$\BGL_{S, \Q} = \oplus_{p \in \Z} \BGL^{(p)}_S$$
such that the pieces $\BGL^{(p)}_S$ represent the graded pieces of the $\gamma$-filtration on $K$-theory:
\eqn
\mylabel{eqn_motivic_vs_K}
\Hom_{\SH(S)}(S^n \wedge \SiginfP X_+, \BGL_S^{(p)}) \cong \gr_\gamma^p K_n(X)_\Q.
\xeqn
%
%
The \Def{Beilinson motivic cohomology spectrum} $\HB$ is defined by
\eqn \mylabel{eqn_motivicCohoSpectrum}
\HBS := \BGL_S^{(0)}
\xeqn
and the resulting Chern character map $\BGL_{S, \Q} \r \oplus_p \HBS\twi p$ is denoted $\ch$.
The parts of the $K$-theory spectrum are related by periodicity isomorphisms
\eqn \mylabel{eqn_BGLtwist}
\BGL^{(p)}_S = \HBS \twi{p}.
\xeqn
For any map $f: T \r S$, not necessarily of finite type, there are natural isomorphisms
\eqn \mylabel{eqn_BGLpullback}
f^* \BGL_S = \BGL_T, \ f^* \HBS = \HBei{T}.
\xeqn

\mylabel{strictringspectra} The following definition and facts are due to Cisinski and D\'eglise \cite[Sections 12.3, 13.2]{CisinskiDeglise:Triangulated}. By a result of R\"ondigs, Spitzweck and Ostvaer \cite{MR2661551}, $\BGL_S \in \SH(S)$ is weakly equivalent to a certain cofibrant strict ring spectrum $\BGL'_S$, that is to say a monoid object in the underlying model category $\Spt^{\P}(\PSh_\bullet(\Sm/S))$. In the same vein, $\HBS$ can be represented by a strict commutative monoid object $\HBS'$ \cite[Cor. 14.2.6]{CisinskiDeglise:Triangulated}. The model structures on the subcategory of $\Spt^{\P}$ of $\BGL'_S$- and $\HBS'$-modules are endowed with model structures such that the forgetful functor is Quillen right adjoint to smashing with $\BGL'_S$ and $\HBS'$, respectively. The homotopy categories are denoted $\DMBGL(S)$ and $\DMBei(S)$, respectively. Objects in $\DMBei(S)$ will be referred to as \emph{motives} over $S$.  We have adjunctions

\eqn \mylabel{eqn_adjunctionModule}
- \wedge \BGL_S : \SH(S) \leftrightarrows \DMBGL(S) : \text{forget}
\xeqn
\eqn \mylabel{eqn_adjunctionModule2}
- \wedge \HBS : \SH(S)_\Q \leftrightarrows \DMBei(S) : \text{forget}.
\xeqn

\mylabel{DMBeiSH} There is a canonical functor from the localization of $\SH(S)_\Q$ by all $\HB$-acyclic objects $E$ (i.e., those satisfying $E \t \HBS = 0$) to $\DMBei(S)$. This functor is an equivalence of categories, which shows that the above definition is independent of the choice of $\HBS'$. This also has the consequence that the forgetful functor $\DMBei(S) \r \SH(S)_\Q$ is fully faithful \cite[Prop. 14.2.8]{CisinskiDeglise:Triangulated}, which will be used in \refsect{definition}. All this stems from the miraculous fact that the multiplication map $\HB \wedge \HB \r \HB$ is an isomorphism.

\Def{Motivic cohomology}  of any object $M$ in $\SH(S)_\Q$ is defined as
\eqn \mylabel{eqn_motiviccohomology}
\H^n(M, p) := \Hom_{\SH(S)_\Q}(M, \HB(p)[n]) \explain{adjunctionModule2} \Hom_{\DMBei(S)}(M \wedge \HBS, \HBS(p)[n]).
\xeqn

The adjunctions \refeq{adjunctionModule}, \refeq{adjunctionModule2} are morphisms of \emph{motivic categories} \cite[Def. 2.4.45]{CisinskiDeglise:Triangulated}, which means in particular that the functors $f_\sharp$, $f_*$, $f^*$, $f_!$ and $f^!$ of \refeq{adjstar}, \refeq{adjshriek}, \refeq{adjsharp} on $\SH(-)$ can be extended to ones on $\DMBGL(-)$ and $\DMBei(-)$ in a way that is compatible with these adjunctions \cite[13.3.3, 14.2.11]{CisinskiDeglise:Triangulated}. For $\DMBei(S)$ this can be rephrased by saying that these functors preserve the subcategories $\DMBei(-) \subset \SH(-)_\Q$.

For any smooth quasi-projective morphism $f : X \r Y$ of constant relative dimension $n$ and any $M \in \DMBei(Y)$, we have the \Def{relative purity} isomorphism (functorial in $M$ and $f$)
\eqn \mylabel{eqn_relativepurity}
f^! M \cong f^* M\twi{n}.
\xeqn
For example, $f^! \HBei{Y} \cong \HBei{X}\twi{n}$. This is due to Ayoub, see e.g.\ \cite[2.4.21]{CisinskiDeglise:Triangulated}.

For any closed immersion $i: X \r Y$ between two regular schemes $X$ and $Y$ with constant relative codimension $n$, there are \Def{absolute purity}  isomorphisms \cite[13.6.3, 14.4.1]{CisinskiDeglise:Triangulated}
\eqn \mylabel{eqn_absolutepurity}
i^! \HBei{Y} \cong \HBei{X} \twi{-n}, \ \ i^! \BGL_Y \cong \BGL_X.
\xeqn

\defi \mylabel{defi_motive}
Let $f: X \r S$ be any map of finite type. We define the \Def{motive} of $X$ over $S$ to be
$$\M(X) := \M_S(X) := f_! f^! \HBei{S} \in \DMBei(S).$$
\xdefi

\rema \mylabel{rema_alternativedefinition}
In \cite[1.1.34]{CisinskiDeglise:Triangulated} the motive of a smooth scheme $f: X \r S$ is defined as $f_\sharp f^* \HBei{S}$. These two definitions agree up to functorial isomorphism: we can assume that $f$ is of constant relative dimension $d$. By relative purity, the functors $f^!$ and  $f^*\twi{d}$ are isomorphic. Thus their left adjoints, namely $f_!$ and $f_\sharp\twi{-d}$ agree, too. Therefore, $f_! f^! \HBei{S} = f_! f^* \HBei{S}\twi{d} = f_\sharp f^* \HBei{S}$.
\xrema

\defi \mylabel{defi_regularprojectivemap}
A map $f : X \r Y$ of $S$-schemes is a \emph{locally complete intersection (l.c.i.) morphism} if both $X$ and $Y$ are regular and, for simplicity of notation, of constant dimension and if
$$f = p \circ i: X \stackrel i \r X' \stackrel p \r Y$$
where $i$ is a closed immersion and $p$ is smooth. Note that this implies that $X'$ is regular. If there is such a factorization with $p: X' = \P[n]_Y \r Y$ the projection, we call $f$ a \emph{regular projective} map.

We shall write $\dim f := \dim X - \dim Y$ for any map $f: X \r Y$ of finite-dimensional schemes.
\xdefi

\exam \mylabel{exam_lci}
Let $f = p \circ i$ be an l.c.i. morphism. Absolute purity for $i$ \refeq{absolutepurity}, relative purity for $p$, and the periodicity isomorphism $\BGL \cong \BGL \twi 1$ give rise to isomorphisms
$$f^! \HBS \cong f^* \HBS \twi{\dim (f)}, \ \ f^! \BGL_S \cong f^* \BGL_S.$$
\xexam

Let $f: X \r Y$ be a projective regular map. Recall the \emph{trace map} in $\SH(Y)$
\eqn \mylabel{eqn_traceBGL}
\tr_f^\BGL: {f}_* \BGL_{X} = {p}_* {i}_* i^* \BGL_{X'} \stackrel{\text{\refeq{absolutepurity}}}\r {p}_* \BGL_{X'} \r \BGL_{Y}
\xeqn
constructed in \cite[13.7.3]{CisinskiDeglise:Triangulated}. This is not an abuse of notation insofar as $\tr_f^\BGL$ is independent of the choice of the factorization. This is shown by adapting \cite[Lemma 5.11]{Deglise:Around} to the case where all schemes in question are merely regular.

The trace map for $\HB$ is defined as the composition
\eqn \mylabel{eqn_traceBeilinson}
\tr_f^{\Beilinson}: f_* f^* \HBei{Y}\twi{\dim f} \rightarrowtail f_* f^* \BGL_{\Q, Y} \stackrel{\tr_f^\BGL}\lr \BGL_{\Q,Y} \twoheadrightarrow  \HBei{Y}.
\xeqn
In case $f = i$, this is the definition of \cite[Section 14.4]{CisinskiDeglise:Triangulated}. 

Given another regular projective map $g$, the composition $g \circ f$ is also of this type. The trace maps are functorial: the composition
$$f^* g^* \BGL \stackrel{\tr_f^\BGL}\lr f^! g^* \BGL \stackrel{f^! \tr_g^\BGL}\lr f^! g^! \BGL$$
agrees with $\tr_{g \circ f}^\BGL$ and similarly with $\tr_?^\Beilinson$. This can be deduced from the independence of the factorization, cf.\ \cite[Prop. 5.14]{Deglise:Around}.

By construction, for any smooth map $f: Y' \r Y$, the induced map $\Hom(f_\sharp f^* S^0, \tr_f^\BGL[-n]) : K_n(X') \r K_n(Y')$ is the $K$-theoretic pushforward along $f' : X' := X \x_Y Y' \r Y'$ \cite[13.7.3]{CisinskiDeglise:Triangulated}. Similarly, $\Hom(f_\sharp f^* S^0, \tr_f^\Beilinson[-n](p))$ is the pushforward $K_n(X')^{(p)}_\Q \r K_n(Y')^{(p)}_\Q$. Indeed, the pushforward on the Adams graded pieces of $K$-theory is defined as the induced map of the graded homomorphism $f'_*$ on $K$-theory \cite[V.6.4]{FultonLang:Riemann}. The adjoint maps
$$\BGL_X = f^* \BGL_Y \r f^! \BGL_Y, \ \  f_* f^* \BGL_Y \r \BGL_Y$$
will also be denoted $\tr_f^\BGL$ and similarly with $\tr_f^\Beilinson$.

\subsection{The Riemann-Roch theorem} \mylabel{sect_RiemannRoch}

We now turn to a motivic Riemann-Roch theorem, which will imply an arithmetic Riemann-Roch theorem for Arakelov motivic cohomology (\refth{HARR}). It generalizes the statement given by Riou for smooth morphisms \cite[Theorem 6.3.1]{Riou:Algebraic} to regular projective maps. Independently, F. D\'eglise has obtained a similar result \cite{Deglise:RiemannRoch}. Recall the \Def{virtual tangent bundle} of a regular projective map $f = p \circ i: X \stackrel i \r X' \stackrel p \r Y$, $T_f := i^* T_p - C_{X/X'} \in K_0(X)$ (see e.g. \cite[V.7]{FultonLang:Riemann}). Here $T_p := \Omega_{X'/Y}^\vee$ is the tangent bundle of $p$ and $C_{X/X'} := (I / I^2)^\vee$ is the conormal sheaf associated the ideal $I$ defining $i$. As an element of $K_0(X)$, $T_f$ does not depend on the factorization. Its Todd class $\Td(T_f)$ is an element of $\oplus_{p\in \Z} K_0(X)_\Q^{(p)}$ (see e.g. \cite[p. 20]{FultonLang:Riemann} for the general definition of $\Td$, this is applied to the Chern character $\ch: K_0(-) \r \oplus_p K_0(-)^{(p)}_\Q$ \cite[p. 127, 146]{FultonLang:Riemann}). It is regarded as an endomorphism of $\oplus_{p \in \Z} \HBei{X}\twi p$ via the natural identification
$\oplus_{p\in \Z} K_0(X)_\Q^{(p)} = \End_{\DMBGL(X)_\Q}(\oplus_{p \in \Z} \HBei{X}\twi p)$.

\theo (Riemann-Roch) \mylabel{theo_RiemannRoch}
Let $f: X \r Y$ be a regular projective map. The following diagram is a commutative diagram in $\SH(Y)_\Q$ (or, equivalently, in $\DMBei(Y)$):
$$\xymatrix{
f_* f^* \BGL_{\Q, Y} \ar[rr]^{\tr_f^\BGL} \ar[d]^\cong_\ch & &
\BGL_{\Q, Y} \ar[d]^\cong_\ch \\
f_* f^* \Beilinson_Y \ar[r]^{f_* \Td(T_f)} &
f_* f^* \Beilinson_Y \ar[r]^{\tr_f^{\Beilinson}} &
{\Beilinson}_Y.}$$
Here, $\Beilinson_Y$ is a shorthand for $\oplus_{p \in \Z} \HBei{Y}\twi{p}$.
\xtheo
\pf
The statement is easily seen to be stable under composition of regular projective maps so it suffices to treat the cases $f = p : \P[n]_Y \r Y$ and $f = i : X \r \P[n]_Y$ separately. The former case has been shown by Riou, so we can assume $f: X \r Y$ is a closed embedding of regular schemes. The classical Riemann-Roch theorem says that the map
$$K_0(X)_\Q \r \oplus_p K_0(Y)^{(p)}_\Q, x \mapsto \ch f_* (x) - f_* (\Td (T_f) \cup \ch(x))$$
vanishes. Viewing $x$ as an element of $\Hom_{\SH(Y)_\Q}(S^0, f_* f^* \BGL_{\Q, Y})$, this can be rephrased by saying that $x \mapsto \alpha_f \circ x$ is zero, where
$$\alpha_f := \ch_X \circ \tr_f^\BGL - \tr_f^\Beilinson \circ f_* \Td(T_f) \circ f_* f^* \ch_Y \in \Hom(f_* f^* \BGL_{\Q, Y}, \Beilinson_{Y}).$$
To show $\alpha_f = 0$, we first reduce to the case where $f : X \r Y$ has a retraction, that is, a map $p: Y \r X$ such that $p \circ f = \id_X$. Then, we prove the theorem by reducing it to the classical Riemann-Roch theorem.

For the first step, recall the deformation to the normal bundle \cite[IV.5]{FultonLang:Riemann}:
\eqn
\mylabel{eqn_deformation}
\xymatrix{
\emptyset \ar[rr] \ar[dd] &
&
X \ar[rr]^{i_\infty} \ar[dd] &
&
\P_X \ar[r]^{pr} \ar[dd]^F &
X \ar[dd]^f \\
&
X \ar@{=}[ur] \ar[dd]^(.3){f'} &
&
X \ar[ur]^{i_0} \ar[dd]^(.3)f \\
\tilde Y \ar[rr] &
&
\tilde Y + Y' \ar[rr]^(.3){s+g'} &
&
M \ar[r]^\pi  &
Y \\
&
Y' \ar[ur] \ar @/^2pc/ @{-->} [uu]
\ar[urrr]_{g'} &
&
Y \ar[ur]^(.3)g \ar@{=}[urr]
}
\xeqn
We have written $M := \text{Bl}_{X \x \infty} (\P_Y)$ and $Y' := \P[](C_{X/Y} \oplus \cO_X)$, $\tilde Y := \text{Bl}_X Y$ and $Y' + \tilde Y$ for the scheme defined by the sum of the two divisors. All schemes except $Y' + \tilde Y$ are regular, all maps except $\pi$ and $pr$ are closed immersions.
The diagram is commutative and every square in it is cartesian. The map $f'$ has a retraction. We show
$$\alpha_{f'} = 0 \Rightarrow \alpha_f = 0$$
by indicating how to replace each argument in \cite[proof of Theorem II.1.3]{FultonLang:Riemann}, which shows $\alpha_{f'} \circ x = 0 \Rightarrow \alpha_f \circ x = 0$ for any $x$ as above, in a manner that is independent of $x$.

The identity $f_* (x) = f_* i_0^* pr^* (x) = g^* F_* pr^* (x)$ is replaced by the commutativity of the following diagram of maps of ($\BGL$-)motives, where $v := g \circ f = F \circ i_0$:
$$\xymatrix{
F_! F^! \BGL_M \ar[rr]^{\cO_{\P_X} \in K_0(\P_X)} \ar[d]_{\cO_X \in K_0(X)} &
&
\BGL_M \ar[d]^{\cO_Y \in K_0(Y)} \\
v_! v^! \BGL_M \ar[rr]_{\cO_X \in K_0(X)} &
&
g_! g^! \BGL_M
}$$
The maps are given by the indicated structural sheaves in $K_0(?)$, via the identifications of $\Hom$-groups in $\DM_\BGL(Y)$ with $K$-theory. For example, the upper horizontal map is the adjoint map to the inverse of the trace map isomorphism $\tr_F^\BGL: F^* \BGL \r F^! \BGL$, which corresponds via absolute purity to $\cO_{\P_X} \in K_0(\P_X) = \Hom_{\DM_\BGL(Y)}(F_! F^! \BGL, \BGL)$. The composition of the map given by $\cO_{\P_X}$ and $\cO_Y$ is given by their tensor product (viewed as $\cO_M$-modules), that is, $\cO_X$, so the diagram commutes. The same argument applies to $f'_* (x) = g'^* F_* pr^* (x)$.

The projection formula is \cite[Theorem 2.4.50(v)]{CisinskiDeglise:Triangulated}. The divisors $Y$ and $Y' + \tilde Y \subset M$ are linearly equivalent, which implies $g_* (1) = g'_* (1) + s_* (1) \in K_0(M)^{(1)}_\Q$ \cite[IV.(5.11), Prop. V.4.4]{FultonLang:Riemann}. This in turn is equivalent to the agreement of the following two elements of $\Hom(\HBei{M}, \HBei{M}\twi{-1})$:
$$\HBei{M} \stackrel{\text{adj.}}\r g_* g^* \HBei{M} \stackrel{g_! \tr_g^\Beilinson}\lr g_! g^! \HBei{M}\twi{-1} \stackrel{\text{adj.}}\r \HBei{M}\twi{-1}$$
and
\small
$$\HBei{M} \stackrel{\text{adj.}}\r g'_* g'^* \HBei{M} \oplus s_* s^* \HBei{M} \stackrel{g'_! \tr_{g'}^\Beilinson \oplus s_! \tr_{s}^\Beilinson} \lr g'_! g'^! \HBei{M}\twi{-1} \oplus s_! s^! \HBei{M}\twi{-1} \stackrel{\text{adj.}}\r \HBei{M}\twi{-1}.$$
\normalsize
Finally, the identity $s^* F_* pr^* (x) = 0$ is formulated independently of $x$ using again base-change (and using that the motive of the empty scheme is zero). This finishes the first step.


Thus, we can assume that $f$ has a retraction $p: Y \r X$. By \cite[Section 5, esp. 5.3.6, cf.\ the proof of 6.1.3.2]{Riou:Algebraic}, the obvious "evaluation" maps  $\Hom(\BGL_{X, \Q}, \BGL_{X, \Q}) $ injectively to
$$\prod_{i \in \Z, T \in \Sm/X} \Hom_{\Q} \left (\Hom((\P)^{\wedge i} \wedge T_+, \BGL_{X, \Q}), \Hom((\P)^{\wedge i} \wedge T_+, \BGL_{X, \Q}) \right).$$
The outer $\Hom$ denotes $\Q$-linear maps, the inner ones are morphisms in $\SH(X)_\Q$. There is an isomorphism $u: f^* \BGL_{\Q, Y} \r f^! \Beilinson_Y$, for example the Chern class followed by the absolute purity isomorphism (\refex{lci}). Appending $u$ on both sides, we conclude that the evaluation maps $\Hom(f^* \BGL_{Y, \Q}, f^! \Beilinson_Y)$ into
$$\prod_{i, T} \Hom_{\Q} \left (\Hom((\P)^{\wedge i} \wedge T_+, f^* \BGL_{Y, \Q}), \Hom((\P)^{\wedge i} \wedge T_+, f^! \Beilinson_Y) \right).$$
For any $T \in \Sm / X$, consider the following cartesian diagram:
$$\xymatrix{
T \ar[r]^{f_T} \ar[d]^t & U \ar[d] \ar[r]^{p_T} & T \ar[d]^t \\
X \ar[r]^f & Y \ar[r]^p & X.
}$$
Recall that $T \in \SH(X)$ is given by $t_\sharp t^* S^0$. Here $t_\sharp$ is left adjoint to $t^*$, cf.\ \refeq{adjsharp}. Thus, the term simplifies to
$$\prod_{i, T} \Hom_{\Q} \left (\Hom((\P)^{\wedge i}, t^* f^* \BGL_{Y, \Q}), \Hom((\P)^{\wedge i}, t^* f^! \Beilinson_Y) \right).$$
The diagram $X \r Y \r X$ is stable with respect to smooth pullback: $f_T$ is also an embedding of regular schemes, $p_T$ is a retract of $f_T$. Moreover, the trace map $\tr^\BGL_f$ behaves well with respect to smooth pullback, i.e., $t^* \tr^\BGL_f = \tr^\BGL_{f_T}$ and similarly for $\tr^\Beilinson_?$, $\ch_?$ and $\Td(T_?)$. Thus, it is sufficient to consider the case $T = X$. That is, we have to show that $\beta_f$, the image of $\alpha_f$ in
\eqnarr
& & \prod_{i \in \Z} \Hom_{\Q} \left (\Hom((\P)^{\wedge i}, f^* \BGL_{Y, \Q}), \Hom((\P)^{\wedge i}, f^! \Beilinson_{Y}) \right) \\
& = & \prod_{i \in \Z} \Hom_{\Q} \left (\Hom_{\SH(X)_\Q}((\P_X)^{\wedge i}, \BGL_{X, \Q}), \Hom_{\SH(Y)_\Q}((\P_Y)^{\wedge i}, f_* f^! \Beilinson_{Y}) \right)
\xeqnarr
is zero. 
The composition
$$\Hom((\P_Y)^{\wedge i}, f_! f^* \Beilinson_Y) \stackrel{\tr_f^\Beilinson, \cong}\lr \Hom((\P_Y)^{\wedge i}, f_! f^! \Beilinson_Y) \stackrel{\gamma_f}\lr \Hom ((\P_Y)^{\wedge i}, \Beilinson_Y)$$
is the pushforward $f_*: \oplus_{p \in \Z} K_0(X)^{(p)}_\Q \r \oplus K_0(Y)^{(p)}_\Q$, which is injective since $p_* f_* = \id$. Thus, the right hand adjunction map $\gamma_f$ is also injective and it is sufficient to show $\gamma_f \circ \beta_f = 0$. For any $i \in \Z$,
\small
\eqnarr
\gamma_f \circ \beta_f & \stackrel{\text{by def.}}= & \left (f_* \circ (- \cup \Td(T_f)) \circ \ch_X \right) - \left (\ch_Y \circ f_* \right) \\
& \stackrel{\text{RR}} = &  0 \\
& \in & \Hom_\Q \left (K_0(X)_\Q, \oplus K_0(Y)^{(p)}_\Q \right) \\
& = & \Hom_{\Q} \left (\Hom_{\SH(X)_\Q}((\P)^{\wedge i}, f^* \BGL_{Y, \Q}), \Hom_{\SH(Y)_\Q}((\P)^{\wedge i}, \Beilinson_{Y}) \right).
\xeqnarr
\normalsize
The vanishing labeled RR is the classical Riemann-Roch theorem for $f$.
\xpf

\subsection{Deligne cohomology} \mylabel{sect_DeligneCohomology}
\defi \cite[3.1.1.]{GilletSoule:Arithmetic} \mylabel{defi_arithmring}
An \Def{arithmetic ring} is a datum $(S, \Sigma, \Fr_\infty)$, where $S$ is a ring, $\Sigma = \{\sigma_1, \dots, \sigma_n : S \r \CC \}$ is a set of embeddings of $S$ into $\CC$ and $\Fr_\infty: \CC^\Sigma \r \CC^\Sigma$ is a $\CC$-antilinear involution (called \emph{infinite Frobenius}) such that $\Fr_\infty \circ \sigma = \sigma$, where $\sigma = (\sigma_i)_i : S \r \CC^\Sigma$. For simplicity, we suppose that $S_\eta  := S \x_{\Spec \Z} \Spec \Q$ is a field. If $S$ happens to be a field itself, we refer to it as an \Def{arithmetic field}. 
For any scheme $X$ over an arithmetic ring $S$, we write
$$X_\CC  := X \x_{S, \sigma} \CC^\Sigma$$
and $X(\CC)$ for the associated complex-analytic space (with its classical topology). We also write $\Fr_\infty: X_\CC \r X_\CC$ for the pullback of infinite Frobenius on the base.
\xdefi

The examples to have in mind are the spectra of number rings, number fields, $\RR$ or $\CC$, equipped with the usual finite set $\Sigma$ of complex embeddings and $\Fr_\infty: (z_v)_{v \in \Sigma} \mapsto (\ol {z_{\ol v}})_v$.

We recall the properties of Deligne cohomology that we need in the sequel. In order to construct a spectrum representing Deligne cohomology in \refsect{Deligne} we recall Burgos' explicit complex whose cohomology groups identify with Deligne cohomology. In the remainder of this subsection, $X / S$ is a smooth scheme (of finite type) over an arithmetic field.

\defi \mylabel{defi_Burgos} \cite[Def. 1.2, Thm. 2.6]{Burgos:Arithmetic}
Let $E^*(X(\CC))$ be the following complex:
\eqn
\mylabel{eqn_differentialforms}
E^*(X(\CC)) := \varinjlim E^*_{\ol X(\CC)}(\log D(\CC)),
\xeqn
where the colimit is over the (directed) category of smooth compactifications $\ol X$ of $X$ such that $D := \ol X \backslash X$ is a divisor with normal crossings. The complex $E^*_{\ol X(\CC)}(\log D(\CC))$ is the complex of $C^\infty$-differential forms on $\ol X(\CC)$ that have at most logarithmic poles along the divisor (see \lcs for details). 
We write $E^*(X) \subset E^*(X(\CC))$ for the subcomplex of elements fixed under the $\ol{\Fr_\infty^*}$-action.
Forms in $E^*(X)$ that are fixed under complex conjugation are referred to as real forms and denoted $E^*_{\RR}(X)$. As usual, a twist is written as $E^*_\RR(X)(p) := (2 \pi i)^p E^*_\RR(X) \subset E^*(X)$. The complex $E^*(X)$ is filtered by
$$F^p E^*(X) := \oplus_{a \geq p, a + b = *} E^{a, b}(X).$$
Let $\Deligne^*(X, p)$ be the complex defined by
$$\Deligne^n(X, p) := \left \{ \begin{array}{ll}
E^{2p+n-1}_\RR(X)(p-1) \cap \oplus_{a + b = 2p+n-1, a, b < p} E^{a, b}(X) & n < 0 \\
E^{2p+n}_\RR(X)(p) \cap \oplus_{a + b = 2p+n, a, b \geq p} E^{a, b}(X) & n \geq 0
\end{array}
\right.
$$
The differential $d_{\Deligne}(x)$, $x \in \Deligne^n(X, p)$, is defined as $-\mathrm{proj} (d x)$ ($n < -1$), $-2\del \ol \del x$ ($n = -1$), and $dx$ ($n \geq 0$). Here $d$ is the standard exterior derivative, and $\mathrm{proj}$ denotes the projection onto the space of forms of the appropriate bidegrees.
We also set
$$\Deligne := \bigoplus_{p \in \Z} \Deligne(p).$$
The pullback of differential forms turns $\Deligne$ into complexes of presheaves on $\Sm / S$. \emph{Deligne cohomology} (with real coefficients) of $X$ is defined as
$$\HD^n(X, p) := \H^{n-2p}(\Deligne(p)(X)).$$
For a scheme $X$ over an arithmetic ring, such that $X_\eta := X \x_S S_\eta$ is smooth (possibly empty), we set $\HD^n(X, p) := \HD^n(X_\eta)$.
\xdefi

Recall that a complex of presheaves $X \mapsto F_*(X)$ on $\Sm / S$ is said to have \emph{\'etale descent} if for any $X \in \Sm / S$ and any \'etale cover $f: Y \r X$ the canonical map
$$F_*(X) \r \Tot (F_*(\dots \r Y \x_X Y \r Y))$$
is a quasi-isomorphism. 
The right hand side is the total complex defined by means of the direct product. (Below we apply it to $F_*(X) = \Deligne(p)(X)$ which is a complex bounded by the dimension of $X$, so that it agrees with the total complex defined using the direct sum in this case.) The total complex is applied to the \v{C}ech nerve. 
At least if $F$ is a complex of presheaves of $\Q$-vector spaces, this is equivalent to the requirement that
$$F_*(X) \r \Tot (F_*(\mathcal Y))$$
is a quasi-isomorphism for any \'etale hypercover $\mathcal Y \r X$. Indeed the latter is equivalent to $F_*$ satisfying Galois descent (as in \refeq{Galoisdescent}) and Nisnevich descent in the sense of hypercovers. The latter is equivalent to the one in the sense of \v{C}ech nerves by the Morel-Voevodsky criterion (see e.g. \cite[Theorem 3.3.2]{CisinskiDeglise:Triangulated}).

\theo \mylabel{theo_DeligneCohomology}
\begin{enumerate}[(i)]
\item \mylabel{item_classicalvsBurgos}
The previous definition of Deligne cohomology agrees with the classical one (for which see e.g. \cite{EsnaultViehweg:Deligne}). In particular, there is a long exact sequence
\eqn \mylabel{eqn_longExactDeligneBettideRham}
\H^n_\Deligne(X, p) \r \H^n(X(\CC), \RR(p))^{(-1)^p} \r (\H^n_\dR(X_\CC) / F^p \H^n_\dR(X_\CC))^{\Fr_\infty} \r \H^{n+1}_\Deligne(X, p)
\xeqn
involving Deligne cohomology, the $(-1)^p$-eigenspace of the $\ol{\Fr_\infty^*}$ action on Betti cohomology, and the $\Fr_\infty$-invariant subspace of de Rham cohomology modulo the Hodge filtration.
\item \mylabel{item_A1invariance}
The complex $\Deligne(p)$ is \emph{homotopy invariant} in the sense that the projection map $X \x \A \r X$ induces a quasi-isomorphism $\Deligne(\A \x X) \r \Deligne(X)$ for any $X \in \Sm /S$.
\item \mylabel{item_firstChern}
There is a functorial \emph{first Chern class} map
\eqn \mylabel{eqn_firstChern}
c_1 : \Pic (X) \r \H^2_\Deligne (X, 1).
\xeqn
\item \mylabel{item_DGA}
The complex $\Deligne$ is a unital differential bigraded $\Q$-algebra which is associative and commutative up to homotopy. The product of two sections will be denoted by $a \cdot_\Deligne b$. The induced product on Deligne cohomology agrees with the classical product $\cup$ on these groups \cite[Section 3]{EsnaultViehweg:Deligne}.
Moreover, for a section $x \in \Deligne^0(X)$ satisfying $d_\Deligne(x) (=dx) = 0$ and any two sections $y, z \in \Deligne^*(X)$, we have
\eqn \mylabel{eqn_productassociative}
x \cdot_\Deligne (y \cdot_\Deligne z) = (x \cdot_\Deligne y) \cdot_\Deligne z
\xeqn
and
\eqn \mylabel{eqn_productcommutative}
x \cdot_\Deligne y = y \cdot_\Deligne x.
\xeqn
\item \mylabel{item_projectiveBundle}
Let $E$ be a vector bundle of rank $r$ over $X$. Let $p: P := \mathbf P(E) \r X$ be the projectivization of $E$ with tautological bundle $\mathcal O_P(-1)$. 
Then there is an isomorphism
\eqn \mylabel{eqn_projectiveBundleDeligne}
p^*(-) \cup c_1 (\mathcal O_P(1))^{\cup i}: \oplus_{i=0}^{r-1} \H^{n-2i}_\Deligne(X, p-i) \r \H^n_\Deligne(P, p).
\xeqn
In particular the following K\"unneth-type formula holds:
\eqn \mylabel{eqn_Kunneth}
\H^n_\Deligne(\P \x X, p) \cong \H^{n-2}_\Deligne(X, p-1) \oplus \H^n_\Deligne(X, p).
\xeqn
\item \mylabel{item_etaleDescent}
The complex of presheaves $\Deligne(p)$ satisfies \'etale descent.
\end{enumerate}
\xtheo
\pf
\refit{classicalvsBurgos}: This explicit presentation of Deligne cohomology is due to Burgos \cite[Prop. 1.3.]{Burgos:Arithmetic}. 
The sequence \refeq{longExactDeligneBettideRham} is a consequence of this and the degeneration of the Hodge to  de Rham spectral sequence. See e.g. \cite[Cor. 2.10]{EsnaultViehweg:Deligne}.
\refit{A1invariance} follows from \refeq{longExactDeligneBettideRham} and the homotopy invariance of Betti cohomology, de Rham cohomology, and, by functoriality of the Hodge filtration, homotopy invariance of $F^p \H^n_\dR (-)$. For \refit{firstChern}, see \cite[Section 5.1.]{BurgosKarmerKuhn} (or \cite[Section 7]{EsnaultViehweg:Deligne} for the case of a proper variety). 
\refit{DGA} is \cite[Theorem 3.3.]{Burgos:Arithmetic}.\footnote{Actually, the product on $\Deligne(X)$ is commutative on the nose. We shall only use the commutativity in the case stated in \refeq{productcommutative} and the associativity as in \refeq{productassociative}, cf.\ \refdele{Delignespectrum}.} 
For \refit{projectiveBundle}, see e.g. \cite[Prop. 8.5.]{EsnaultViehweg:Deligne}.

\refit{etaleDescent}: This statement can be read off the existence of the absolute Hodge realization functor \cite[Corollary 2.3.5]{Huber:Realization} (and also seems to be folklore). Since it is crucial for us in \refth{modulestructureDeligne}, we give a proof here. Let
$$\tilde \Deligne^*(X, p) := \cone (E_\RR^*(X)(p) \oplus F^p E^*(X) \stackrel{(+1, -1)}\lr E^*(X))[-1+2p].$$ 
By \cite[Theorem 2.6.]{Burgos:Arithmetic}, there is a natural (fairly concrete) homotopy equivalence between the complexes of presheaves $\tilde \Deligne(p)$ and $\Deligne(p)$.
The descent statement is stable under quasi-isomorphisms of complexes of presheaves and cones of maps of such complexes. Therefore it is sufficient to
show descent for the complexes $E_\RR^*(-)(p)$, $F^p E^*(-)$, $E^*(-)$.
Taking invariants of these complex under the $\ol {\Fr^*_\infty}$-action is an exact functor, so we can disregard that operation in the sequel. From now on, everything refers to the analytic topology, in particular we just write $X$ for $X(\CC)$ etc.
Let $j: X \r \ol X$ be an open immersion into a smooth compactification such that $D := \ol {X} \backslash X$ is a divisor with normal crossings. The inclusion
$$\Omega^*_{\ol {X}}(\log D) \subset E^*_X(\log D)$$
of holomorphic forms into $C^\infty$-forms (both with logarithmic poles) yields quasi-isomorphisms of complexes of vector spaces
$$\RG \R j_* \CC \r \RG \R j_* \Omega^*_{X} \leftarrow  \RG \Omega^*_{\ol {X}}(\log D) \r \Gamma E^*_{X}(\log D)$$
that are compatible with both the real structure and the Hodge filtration \cite[Theorem 2.1.]{Burgos:Dolbeault}, \cite[3.1.7, 3.1.8]{Deligne:Hodge2}. Here $(\R) \Gamma$ denotes the (total derived functor of the) global section functor on $\ol {X}$.
The complex $E^*(X)$, whose cohomology is $\H^*(X, \CC)$, is known to satisfy \'etale descent \cite[Prop. 2.1.7]{Huber:Realization}. This also applies to $E^*_\RR(X)(p)$ instead of $E^*(X)$. 
(Alternatively for the former, see also \cite[3.1.3]{CisinskiDeglise:Mixed} for the \'etale descent of the algebraic de Rham complex $\Omega^*_X$.)

It remains to show the descent for $X \mapsto F^p E^*(X)$. Consider a distinguished square in $\Sm / S$,
$$\xymatrix{
X' \ar[r] \ar[d] & X \ar[d] \\
Y' \ar[r] & Y,}
$$
i.e., cartesian such that $Y' \r Y$ is an open immersion, $X/Y$ is \'etale and induces an isomorphism $(X \backslash X')_\red \r (Y \backslash Y')_\red$. Then the sequence
\eqn
\mylabel{eqn_longHodge}
\H^n(F^p E^*(Y)) \r \H^n(F^p E^*(Y')) \oplus \H^n(F^p E^*(X)) \r \H^n(F^p E^*(X')) \r \H^{n+1}(F^p E^*(Y))
\xeqn
is exact: first, the direct limit in \refeq{differentialforms} is exact. Moreover, $\H^n(\Gamma(F^p E_{\ol X}(\log D)))$ maps injectively into $\H^n(\ol X, \Omega^*_{\ol X}(\log D))$, and the image is precisely the $p$-th filtration step of the Hodge filtration on $\H^n(\ol X, \Omega^*_{\ol X}(\log D)) = \H^n(X, \CC)$. Similarly for $X'$ etc., so that the exactness of \refeq{longHodge} results from the one of the sequence featuring the Betti cohomology groups of $Y$, $Y' \sqcup X$ and $X'$, together with the strictness of the Hodge filtration  \cite[Th.\ 1.2.10]{Deligne:Hodge2}. This shows Nisnevich descent for the Hodge filtration. Secondly, for any scheme $X$ and a Galois cover $Y \r X$ with group $G$, the pullback map into the $G$-invariant subspace
\eqn \mylabel{eqn_Galoisdescent}
\H^n(F^p E^*(X)) \r \H^n(F^p E^*(Y)^G)
\xeqn
is an isomorphism. Indeed, the similar statement holds for $E^*(-)$ instead of $F^p E^*(-)$. We work with $\Q$-coefficients, so taking $G$-invariants is an exact functor, hence $\H^n(F^p E^*(Y)^G) = (\H^n(F^p E^*(Y)))^G = (F^p \H^n_\dR(Y))^G = F^p (\H^n_\dR(Y)^G)$, the last equality by functoriality of the Hodge filtration. Then, again using the strictness of the Hodge filtration, the claim follows. Hence the presheaf $X \mapsto F^p E^*(X)$ has \'etale descent.
\xpf


\section{The Deligne cohomology spectrum} \mylabel{sect_Deligne}
Let $S$ be a smooth scheme (of finite type) over an arithmetic field (\refde{arithmring}). The aim of this section is to construct a ring spectrum in $\SH(S)$ which represents Deligne cohomology for smooth schemes $X$ over $S$.
The method is a slight variation of the method of Cisinski and Deglise used in \cite{CisinskiDeglise:Mixed} to construct a spectrum for any mixed Weil cohomology, such as algebraic or analytic de Rham cohomology, Betti cohomology, and (geometric) \'etale cohomology. The difference compared to their setting is that the Tate twist on Deligne cohomology groups is not an isomorphism of vector spaces.

In this section, all complexes of (presheaves of) abelian groups are considered with homological indexing: the degree of the differential is $-1$ and  $C[1]$ is the complex whose $n$-th group is $C_{n+1}$. 
As usual, any cohomological complex is understood as a homological one by relabeling the indices.
In particular, we apply this to (the restriction to $\Sm/S$ of) the complexes $\Deligne(p)$, $\Deligne$ defined in \ref{defi_Burgos} and let
\eqn
\mylabel{eqn_Delignehomological}
\Deligne_n := \Deligne^{-n} = \oplus_{p \in \Z} \Deligne^{-n}(p).
\xeqn
In order to have a complex of simplicial presheaves (as opposed to a complex of abelian groups), we use the Dold-Kan-equivalence
$$\cK: \Com_{\geq 0}(\Ab) \rightleftarrows \Delta^\op(\Ab) : \mathcal N$$
between homological complexes concentrated in degrees $\geq 0$ and simplicial abelian groups.
We write $\tau_{\geq n}$ for the good truncation of a complex.

\defi \mylabel{defi_Delignesimplicial}
We write
$$\Deligne_s := \cK (\tau_{\geq 0} \Deligne),$$
$$\Deligne_s(p) := \cK (\tau_{\geq 0} \Deligne(p)).$$
\xdefi

Via the Alexander-Whitney map, the product on $\Deligne$ transfers to a map
\eqn \mylabel{eqn_wedgeDeligne}
\Deligne_s(p) \wedge \Deligne_s(p') \r \Deligne_s(p+p').
\xeqn

\lemm \mylabel{lemm_DeligneHoC}
For $X$ smooth over $S$ and any $k \geq 0$, $p \in \Z$ we have:
\eqn
\mylabel{eqn_DeligneHoC}
\Hom_{\HoC_\bullet}(S^k \wedge X_+, \Deligne_s(p)) = \H^{2p-k}_\Deligne(X, p)
\xeqn
and similarly for $\Deligne_s$.
\xlemm
\pf
In $\HoCsect$ (cf.\ \refsect{SH} for the notation), the $\Hom$-group reads
\begin{eqnarray*}
\Hom_{\HoCsect}(S^k \wedge X_+, \cK (\tau_{\geq 0} (\Deligne))) & = & \pi_{k} \cK (\tau_{\geq 0} (\Deligne(X))) \\
& = & \H_{k} (\tau_{\geq 0} (\Deligne(X))) \\
& = & \oplus_{p \in \Z} \H^{2p-k}_\Deligne (X, p).
\end{eqnarray*}
%
%
We have used the fact that any simplicial abelian group is a fibrant simplicial set and the identification $\pi_n (A, 0) = \H_n(\mathcal N(A))$ for any simplicial abelian group.

The presheaf $\Deligne_s$ is fibrant with respect to the motivic model structure, since Deligne cohomology satisfies Nisnevich descent and is $\A$-invariant by
\refth{DeligneCohomology} \refit{etaleDescent} and \refit{A1invariance}. Thus the $\Hom$-groups agree when taken in $\HoCsect$ and $\HoC$, respectively.
\xpf

\defilemm \mylabel{defilemm_Delignespectrum}
The \Def{Deligne cohomology spectrum} $\cDz$ is the spectrum consisting of the $\Deligne_s(p)$ ($p \geq 0$), equipped with the trivial action of the symmetric group $\Sigma_p$. We define the bonding maps to be the composition
$$\sigma_p : \P_S \wedge \Deligne_s(p) \stackrel{c^* \wedge \id} \r \Deligne_s(1) \wedge \Deligne_s(p) \stackrel{\text{\refeq{wedgeDeligne}}}{\to} \Deligne_s(p+1).
$$
Here $c^*$ is the map induced by $c := c_1(\cO_{\P}(1), FS) \in \Deligne_0(1)(\P)$, the first Chern form of the bundle $\cO(1)$ equipped with the Fubini-Study metric. This defines a symmetric $\P$-spectrum.

Define the unit map $1_{\Deligne}: \SiginfP S_+ \to \cDz$ in degree zero by the 
unit of the DGA $\Deligne(0)$. 
In higher degrees, we put
\eqn \mylabel{eqn_UnitDeligne}
(1_{\Deligne})_p: (\P)^{\wedge p} \stackrel{(c^*)^{\wedge p}} \lr \Deligne_s(1)^{\wedge p} \stackrel{\mu}\lr \Deligne_s(p).
\xeqn
Equivalently, $(1_{\Deligne})_p := \sigma_{p-1} \circ (\id_{\P} \wedge (1_{\Deligne})_{p-1})$. This map and the product map $\mu: \cDz \wedge \cDz \r \cDz$ induced by \refeq{wedgeDeligne} turns $\cDz$ into a commutative monoid object of $\SH(S)$, i.e., a commutative ring spectrum. 
\xdefilemm
\pf
Recall that $c$ is a $(1,1)$-form which is invariant under $\Fr_\infty^*$ and under complex conjugation, so $c$ is indeed an element of $\Deligne_0(1)(\P)$. Its restriction to the point $\infty$ is zero for dimension reasons, so $c$ is a pointed map $(\P, \infty) \r (\Deligne_0(1), 0)$. It remains to show that the map
\eqnarr
(\P)^{\wedge m} \wedge \Deligne_s (n) & \stackrel{\id^{\wedge m-1} \wedge c^* \wedge \id} \lr &
(\P)^{\wedge m-1} \wedge \Deligne_s(1) \wedge \Deligne_s (n)  \\
& \stackrel{\text{\refeq{wedgeDeligne}}}\r & (\P)^{\wedge m-1} \wedge \Deligne_s(n+1) \\
& \r & \dots \\
& \r & \Deligne_s(m+n)
\xeqnarr
is $\Sigma_m \x \Sigma_n$-equivariant, i.e., invariant under permuting the $m$ wedge factors $\P$. Given some map $f: U \r (\P)^{\x m}$ with $U \in \Sm/S$, let $f_i: U \r \P$ be the $i$-th projection of $f$ and $c_i := f_i^* c_1(\cO_{\P}(1))$. Given some form $\omega \in \Deligne(n)(U)_*$, we have to check that the expression
\eqn 
\mylabel{eqn_c1comm}
c_1 \cdot_\Deligne (c_2 \cdot_\Deligne (\dots (c_m \cdot_\Deligne \omega)  \dots ))
\xeqn
is invariant under permutation of the $c_i$. Here $\cdot_\Deligne$ stands for the product map \refeq{wedgeDeligne}. This holds before applying the Dold-Kan functor $\mathcal K$ (i.e., $(\P)^{\x m} \x \Deligne(n) \r \Deligne(n+m)$ is $\Sigma_m$-invariant) since the forms $c_i \in \Deligne_0(1)(U)$ are closed, so by \refth{DeligneCohomology}\refit{DGA} the expression \refeq{c1comm} is associative and commutative. The Alexander-Whitney map is symmetric in (simplicial) degree $0$, i.e. $\mathcal K (\Deligne(p)) \wedge \mathcal K(\Deligne(p')) \r \mathcal K (\Deligne(p) \t \Deligne(p'))$ commutes with the permutation of the two factors when restricting to elements of degree $0$. Moreover, it is associative in all degrees. As $c_i \in \Deligne_0(1)$, the previous argument carries over to the product on $\Deligne_s(-)$ instead of $\Deligne(-)$. This shows that $\HD$ is a symmetric spectrum.   

By \lc, the product on $\Deligne$ is (graded) commutative and associative up to homotopy, thus the diagrams checking, say, the commutativity of $\cDz \wedge \cDz \r \cDz$ do hold in $\SH(S)$. The details of that verification are omitted.
\xpf


\rema 
\begin{enumerate}
\item Consider the spectrum $D'$ obtained in the same way as $\HD$, but replacing $\Deligne_s(p)$ by $\cDz$. Then the obvious map $\phi: \bigoplus_{p \in \Z} \cDz\twi{p} \r D'$ is an isomorphism. To see that, it is enough to check that $\Hom_{\SH(S)}(S^n \wedge \SiginfP X_+ , -)$ yields an isomorphism when applied to $\phi$. By the compactness of $S^n \wedge \SiginfP X_+$ in $\SH(S)$, this $\Hom$-group commutes with the direct sum. Then the claim is trivial.
\item
Choosing another metric $\lambda$ on $\cO(1)$ in the above definition, the resulting Deligne cohomology spectrum would be weakly equivalent to $\cDz$ since the difference of the Chern forms $c_1(\cO(1), FS) - c_1(\cO(1), \lambda)$ lies in the image of $d_\Deligne : D_1(1) \r D_0(1)$, see e.g. \cite[Lemma 5.6.1]{Jost:Compact}.
\end{enumerate}
\xrema

\lemm \mylabel{lemm_DeligneOmega}
The Deligne cohomology spectrum $\cDz$ is an $\Omega$-spectrum (with respect to smashing with $\P$). 
\xlemm
\pf
We have to check that the adjoint map to $\sigma_p$ (\refdele{Delignespectrum}),
$$b_p : \Deligne_s(p) \r \R\IHom_\bullet(\P, \Deligne_s(p+1)),$$
is a motivic weak equivalence. As $\P$ is cofibrant and $\Deligne_s(p+1)$ is fibrant, the
non-derived $\IHom_\bullet(\P, \Deligne_s(p))$ is fibrant and agrees with $\R\IHom_\bullet(\P, \Deligne_s(p))$. The map is actually a sectionwise weak equivalence, i.e., an isomorphism in $\HoCsect(S)$. To see this, it is enough to check that the map
$$\Deligne_s(p)(U) \r \IHom_\bullet (\P, \Deligne_s(p+1)(U))$$
is a weak equivalence of simplicial sets for all $U \in \Sm / S$ \cite[1.8., 1.10, p. 50]{MorelVoevodsky:A1}. The $m$-th homotopy group of the left hand side is $\H^{2p-m}_\Deligne(U, p)$ (\refle{DeligneHoC}), while $\pi_m$ of the right hand simplicial set identifies with those elements of $\pi_m (\IHom(\P \x U, \Deligne_s(p+1))) =  \H^{2(p+1)-m}_\Deligne(\P \x U, p+1)$ which restrict to zero when applying the restriction to the point $\infty \r \P$. By the projective bundle formula \refeq{Kunneth}, the two terms agree.
\xpf

\theo \mylabel{theo_Delignespectrum} \mylabel{theo_modulestructureDeligne}
\begin{enumerate}[(i)]
\item \mylabel{item_proofi}
The ring spectrum $\cDz$ represents Deligne cohomology in $\SH(S)$: for any smooth scheme $X$ over $S$, and any $n$, $m \in \Z$ we have
$$\Hom_{\SH(S)} ( (S^1)^{\wedge n} \wedge (\P_S)^{\wedge m} \wedge  \SiginfP X_+, \cDz)  = \H_{\Deligne}^{-n-2m} (X, -m).$$
(See p.\ \pageref{smashinginvertible} for the meaning of $(S^1)^{\wedge n} $, $(\P_S)^{\wedge m}$ with negative exponents.)
\item \mylabel{item_proofii}
The Deligne cohomology spectrum $\cDz$ has a unique structure of an $\HBS$-algebra and $\oplus_{p \in \Z} \cDz\twi{p}$ has a unique structure of an $\BGL_S$-algebra. In particular, $\cDz$ is an object in $\DMBei(S)$, so that \refit{proofi} and \refeq{adjunctionModule2} yield a natural isomorphism
$$\Hom_{\DMBei(S)} (\M_S(X), \HD (p)[n]) = \HD^n(X, p)$$
for any smooth $X/S$.
\item \mylabel{item_proofiii}
The map $\id_\Deligne \wedge 1_{\HB}: \HD \r \HD \wedge \HB$ is an isomorphism in $\SH(S)_\Q$.
\end{enumerate}
\xtheo

\defi \mylabel{defi_regulators}
The maps induced by the unit of $\cDz$ are denoted $\rho_\Deligne: \HB \r \cDz$ and $\ch_\Deligne: \BGL \r \oplus_p \HD \twi p$, respectively. 
\xdefi

\pf
By \refle{DeligneOmega}, $\cDz$ is an $\Omega$-spectrum. Thus \refit{proofi} follows from \refle{DeligneHoC}.

\refit{proofii}: by \ref{defilemm_Delignespectrum}, $\cDz$ is a commutative ring spectrum. Recall the definition of \'etale descent for spectra and that for this it is sufficient that the individual pieces of the spectrum have \'etale descent \cite[Def. 3.2.5, Cor. 3.2.18]{CisinskiDeglise:Triangulated}. Thus, $\cDz$ satisfies \'etale descent by \refth{DeligneCohomology}\refit{etaleDescent}. Moreover, $\cDz$ is orientable since $\Hom_{\SH(S)}(\P[\infty]_S, \HD\twi 1) = \varprojlim_n \Hom(\P[n], \HD \twi 1)$ by the Milnor short exact sequence (see e.g. \cite[Cor. 2.2.8]{CisinskiDeglise:Mixed} for a similar situation). This term equals $\HD^2(\P, 1)$ by \refeq{projectiveBundleDeligne}. Any object in $\SH(S)_\Q$ satisfying \'etale descent is an object of $\DMBei(S)$, i.e., an $\HBS$-module \cite[proof of 16.2.18]{CisinskiDeglise:Triangulated}. If it is in addition an orientable ring spectrum, there is a unique $\HBS$-algebra structure on it \cite[Cor. 14.2.16]{CisinskiDeglise:Triangulated}. This settles the claim for $\cDz$. Secondly, the natural map  (in $\SH(S)$)
$$\BGL \r \BGL_\Q \stackrel{\text{\refeq{BGLtwist}}} \cong \oplus_{p \in \Z} \HB\twi{p} \stackrel{\rho_\Deligne\twi{p}}\lr \oplus_p \cDz\twi{p}$$
and the ring structure of $\oplus \cDz\twi{p}$ defines a $\BGL$-algebra structure on $\oplus \cDz\twi{p}$. This uses that the isomorphism \refeq{BGLtwist} is an isomorphism of monoid objects \cite[14.2.17]{CisinskiDeglise:Triangulated}. The unicity of that structure follows from the unicity of the one on $\cDz$ and $\Hom_{\SH(S)}(\BGL_\Q, \oplus \cDz\twi{p}) = \Hom_{\SH(S)_\Q}(\BGL_\Q, \oplus \cDz\twi{p})$, since $\cDz$ is a spectrum of $\RR$- (a fortiori: $\Q$\nobreakdash-)vector spaces.

\refit{proofiii} follows from \refit{proofii}, using \cite[14.2.16]{CisinskiDeglise:Triangulated}.
\xpf

\section{Arakelov motivic cohomology} \mylabel{sect_Arakelov}
Let $S$ be a regular scheme of finite type over an arithmetic ring $B$. The generic fiber $S_\eta := S \x_\Z \Q \r B_\eta := B \x_\Z \Q$ is smooth, since $B_\eta$ is a field (by \refde{arithmring}). We now define the Arakelov motivic cohomology spectrum $\HBhatS$ which glues, in a sense, the Deligne cohomology spectrum $\cDz \in \SH(S_\eta)$ (\refsect{Deligne}) with the Beilinson motivic cohomology spectrum $\HBS$ \refeq{motivicCohoSpectrum}. Parallely, we do a similar construction with $\BGL_S$ instead of $\HBS$.
Once this is done, the framework of the stable homotopy category and motives readily imply the existence of functorial pullbacks and pushforwards for Arakelov motivic cohomology (\refsect{functoriality}). We also prove a higher arithmetic Riemann-Roch theorem (\refth{HARR}) and deduce further standard properties, such as the projective bundle formula in \refsect{further}.

\subsection{Definition} \mylabel{sect_definition}
Recall from \refsect{SH} the category $\Spt(S) := \Spt^{\P}(\Delta^\op \PSh_\bullet(\Sm/S))$ with the stable model structure. The resulting homotopy category is $\SH(S)$.

\defi \mylabel{defi_widehatA}
For any $A \in \Spt(S)$, we put
\eqn
\mylabel{eqn_definitionArakelov}
\widehat A := \hofib_{\Spt(S)} \left (A \wedge \mathrm Q \R (S^0) \stackrel{\id \wedge \mathrm Q \R (1_\Deligne)} \lr A \wedge \mathrm{Q}\R \eta_* \HD \right) \in \Spt(S).
\xeqn
Here, $\hofib$ stands for the homotopy fiber, $1_\Deligne : S^0 \r \HD$ is the unit map given in \refeq{UnitDeligne}, and $\mathrm Q$ and $\R$ are the cofibrant and fibrant replacement functors in $\Spt(S)$. The map $1_\Deligne$ is a map in $\Spt(S_\eta)$, as opposed to a map in the homotopy category $\SH(S_\eta)$. Hence so is the map used in \refeq{definitionArakelov}. We wrote $\mathrm Q \R$ here for clarity, but drop these below, given that the fibrant-cofibrant replacement of any spectrum is weakly equivalent to the original one.

We write $[\widehat A]$ for the image of $\widehat A$ in $\SH(S)$ (or $\SH(S)_\Q$) under the localization functor. Using the strict ring spectra $\HBS'$ and $\BGL'_S$ (p.\ \pageref{strictringspectra}), we define the \Def{Arakelov motivic cohomology spectrum} $\HBhatS$ as 
$$\HBhatS := [\widehat {\HBS'}]  \in \SH(S)_\Q$$
and similarly
$$\BGLhatS := [\widehat {\BGL'}] \in \SH(S).$$
\xdefi

\theo \mylabel{theo_hats}
\begin{enumerate}[(i)]
\item \mylabel{item_hats1} Given a morphism $f: A \r A'$ in $\SH$, there is a canonical morphism $[\widehat f]: [\widehat A] \r [\widehat {A'}]$ in $\SH$ which is an isomorphism if $f$ is. In particular, the Chern character isomorphism $\ch: \BGL_{S, \Q} \cong \oplus_{p \in \Z} \HBei{S}\twi p$ gives rise to an isomorphism called \emph{Arakelov Chern character},
\eqn \mylabel{eqn_ArakelovChernCharacter}
\widehat {\ch}: \BGLhat_{S, \Q} \cong \oplus \HBhatS \twi p
\xeqn
in $\SH(S)_\Q$.
\item \mylabel{item_hats2} If $A$ is a strict ring spectrum, then $[\widehat A]$ is an $A$-module in a canonical way. In particular, $\HBhatS$ is in $\DMBei(S)$ and $\BGLhatS$ is an object of $\DMBGL(S)$. 
\end{enumerate}
\xtheo

\pf
\refit{hats1}: we can represent $f$ by a zig-zag of maps $f_i$ and define $\widehat f$ to be the zig-zag of $\widehat {f_i} := \hofib (f_i \wedge (S^0 \stackrel {1_{\HD}} \r \HD))$. As any choice of the zig-zag represents the same given map $[f]: [A] \r [A']$ in $\SH(S)$, the resulting map $[\widehat f] : [\widehat A] \r [\widehat {A'}]$ is also independent of the choice of the zig-zag.  

\refit{hats2}: the map in \refeq{definitionArakelov} is a map of $A$-modules. Its homotopy fiber in the category of $A$-modules is an object $\widehat{A^{\Mod}} \in A-\Mod$. By the Quillen adjunction \refeq{adjunctionModule2} and \cite[Theorem 19.4.5]{Hirschhorn:Model}, $\widehat{A^{\Mod}}$ is weakly equivalent (in $\Spt$) to $\widehat A$. Therefore, the image of $[\widehat{A^{\Mod}}]$ in $\SH$ under the forgetful functor $\Ho(A-\Mod) \r \SH$ is isomorphic to $[\widehat A]$, i.e.,  latter is canonically an $A$-module.
\xpf

\rema \mylabel{rema_propertiesHBhat}
\begin{enumerate}[(i)]
\item \mylabel{item_hats3} \refth{hats}\refit{hats1} shows that $\BGLhat$ does not depend on the choice of the spectrum representing $\BGL$. In a similar vein, one can show that given a map $A \r A'$ of strict ring spectra (respecting the ring structure) that is also a weak equivalence, $[\widehat A]$ is mapped to $[\widehat {A'}]$ under the canonical equivalence of categories $-: \t^L_A A' : \Ho(A-\Mod) \r \Ho(A'-\Mod)$. In this sense, the $\BGL$-module structure on $\BGLhat$ does not depend on the choice of the strict ring spectrum. We will not use this fact, though.
\item 
We are mainly interested in gluing motivic cohomology with Deligne cohomology. However, nothing is special about Deligne cohomology. In fact, given some scheme $f: T \r S$ (not necessarily of finite type), and complexes of presheaves of $\Q$-vector spaces $D(p)$ on $\Sm / T$ satisfying the conclusion of \refth{DeligneCohomology}\refit{A1invariance}, \refit{firstChern}, \refit{projectiveBundle} (actually \refeq{Kunneth} suffices), \refit{etaleDescent}, and \refit{DGA}, everything could be done with $f_* D(p)$ instead of $\eta_* \Deligne(p)$.
\end{enumerate}
\xrema

\defi \mylabel{defi_Hhats}
For any $M \in \SH(S)$, we define
$$\Hhat^n(M) := \Hom_{\SH(S)}(M, \BGLhatS[n]),$$
$$\Hhat^n(M, p) := \Hom_{\SH(S)_\Q}(M_\Q, \HBhat(p)[n]).$$
The latter is called \Def{Arakelov motivic cohomology} of $M$. For any finite type scheme $f : X \r S$, we define Arakelov motivic cohomology of $X$ as
$$\Hhat^n(X/S, p) := \Hom_{\SH(S)_\Q}(f_! f^! \Sigma_{\P}^\infty S^0, \HBhatS(p)[n])$$
and likewise
$$\Hhat^n(X/S) := \Hom_{\SH(S)}(f_! f^! \Sigma_{\P}^\infty S^0, \BGLhatS[n]).$$
Here $\Sigma_{\P}^\infty S^0$ is the infinite $\P$-suspension of the $0$-sphere, i.e., the unit of the monoidal structure in $\SH$. When the base $S$ is clear from the context, we will just write $\Hhat^n(X, p)$ and $\Hhat^n(X)$. See \refth{HARR}\refit{purity} for a statement concerning the independence of the base scheme $S$ of the groups $\Hhat^n(X/S)$.
\xdefi

\theo \mylabel{theo_longexactArakelov}
\begin{enumerate}[(i)]
\item \mylabel{item_proof1_i}
For any $M \in \SH(S)$ there are long exact sequences relating Arakelov motivic cohomology to (usual) motivic cohomology \refeq{motiviccohomology} and, for appropriate motives, Deligne cohomology (\refde{Burgos}):
\eqn \mylabel{eqn_longExactArakelov}
\dots \r \Hhat^n (M, p) \r \H^n (M, p) \stackrel \rho \r \Hom_{\SH(S)}(M, \eta_* \cDz(p)[n]) \r  \Hhat^{n+1} (M, p) \dots
\xeqn
\eqn 
\dots \r \Hhat^n (M) \r \H^n (M) \stackrel \ch \r \Hom_{\SH(S)}(M, \oplus \eta_* \cDz \twi p[n]) \r  \Hhat^{n+1} (M) \dots
\xeqn
The maps $\rho$ and $\ch$ agree with the one induced by $\rho_\Deligne$ and $\ch_\Deligne$ (\refde{regulators}).
\item \mylabel{item_proof1_ii}
For any l.c.i.\ scheme $X/S$ (\refde{regularprojectivemap}, for example $X = S$) we get exact sequences
$$\cdots \r \Hhat^n (X, p) \to K_{2p-n}(X)^{(p)}_\Q \to \HD^n (X, p) \to  \Hhat^{n+1} (X, p) \r \cdots,$$
\eqn \mylabel{eqn_longexactK}
\cdots \r \Hhat^n(X) \r K_{-n}(X) \r \oplus_p \HD^{2p-n} (X, p) \r \Hhat^{n+1}(X) \r \cdots .
\xeqn
\item \mylabel{item_proof1_iii}
If $S' \stackrel f \r S$ is a scheme of positive characteristic over $S$, the obvious map $\Hhat^n (f_* M, p) \r \H^n (f_* M, p)$ is an isomorphism for any $M \in \SH(S')$.
\item \mylabel{item_proof1_iv}
There is a functorial isomorphism 
\eqn \mylabel{eqn_HhatModule}
\Hhat^n(M) = \Hom_{\DMBGL(S)} (\BGL_S \wedge M, \BGLhatS),
\xeqn
where we view $\BGLhatS$ as a $\BGL$-module using \refth{hats}\refit{hats2}. A similar statement holds for $\HBS$. In addition, there is a canonical isomorphism $\Hhat^n(M, p) = \Hhat^n(M \wedge \HBS, p)$. For example, $\Hhat^n(X, p) = \Hhat^n(\M_S(X), p)$ for any $X/S$ of finite type. For any compact object $M \in \SH(S)$, there is an isomorphism called the \emph{Arakelov Chern character}:
\eqn \mylabel{eqn_HBhatVsBGLhat}
\widehat{{\ch}}: \Hhat^n(M) \t_\Z \Q = \oplus_{p \in \Z} \Hhat^{n+2p}(M, p).
\xeqn
\end{enumerate}
\xtheo
\pf
The long exact sequence in \refit{proof1_i} follows from \refth{Delignespectrum}\refit{proofiii}, the projection formula $\HB \wedge \eta_* \HD = \eta_* (\HB \wedge \HD)$, and generalities on the homotopy fiber in stable model categories.  
Similarly, $\BGL \wedge \HD$ is canonically isomorphic, via the Chern class, to $\oplus \HB \wedge \HD \twi p \cong \oplus \HD \twi p$. The agreement of $\rho$ and $\rho_\Deligne$ is also clear, since the $\HB$-module structure map $\HB \wedge \HD \r \HD$ is inverse to $1_\Beilinson \wedge \id_\Deligne : \HD \r \HB \wedge \HD$.

For \refit{proof1_ii}, we use \refit{proof1_iv} and apply \refit{proof1_i} to $\M_S(X)$ and $f_! f^! \BGL_S$, respectively, where $f: X \r S$ is the structural map. In order to identify the motivic cohomology with the Adams eigenspace in $K$-theory, we use the adjunction \refeq{adjshriek} and the purity isomorphism for $f$ (\refex{lci}). 
To calculate $\Hom(f_! f^! \HBS, \HD)$, we can replace $B$ by the arithmetic field $B_\eta := B \x_\Z \Q$. 
The scheme $S$ is regular, thus $s: S \r B$ is smooth (of finite type). The same is true for the structural map $x: X \r B$. Now, combining the relative purity isomorphisms for $x$ and for $s$, we get an isomorphism
\eqnarr
f^! \HD & = & f^! s^* \HD = f^! s^! \HD \twi{-\dim s} \\
& = & x^! \HD \twi {-\dim s} = x^* \HD \twi {-\dim s + \dim x} = f^* \HD \twi{\dim f}.
\xeqnarr
We conclude
\eqnarr
\Hom_{\SH(S)} (f_! f^! \HBS, \cDz(p)[n]) & = & \Hom(f^! \HBei{S}, f^! \cDz(p)[n]) \\
& = & \Hom(f^* \HBei{S} \twi{\dim f}, f^* \cDz(p)[n] \twi{\dim f}) \\
& = & \Hom(\HBei{X}, \cDz(p)[n]) \\
& \stackrel{\text{\ref{theo_Delignespectrum}}}= & \HD^{2p-n}(X, p).
\xeqnarr
\refit{proof1_iii} follows from localization. The first isomorphism in \refit{proof1_iv} follows from \refeq{adjunctionModule2}. The second one uses in addition the full faithfulness of the forgetful functor $\DMBei \r \SH_\Q$ (p.\ \pageref{DMBeiSH}). The map $\widehat{\ch}$ is induced by \refeq{ArakelovChernCharacter}.
\xpf

\rema
By \refeq{longExactArakelov}, each group $\Hhat^n(M)$ is an extension of a $\Z$-module by a quotient of a finite-dimensional $\RR$-vector space by some $\Z$-module. Both $\Z$-modules are conjectured to be finitely generated in case $S = \SpecZ$ and $M$ compact (Bass conjecture). Similarly, the groups $\Hhat^n(M, p)$ are extensions of $\Q$-vector spaces by groups of the form $\RR^k / \text{some }\Q\text{-subspace}$. In particular, we note that the Arakelov motivic cohomology groups $\Hhat^n(M, p)$ are typically infinite-dimensional (as $\Q$-vector spaces). However, one can can redo the above construction using the spectrum $\HB \otimes \RR$ instead of $\HB$ to obtain \emph{Arakelov motivic cohomology groups with real coefficients}, $\Hhat^n(M, \RR(p))$. These groups are real vector spaces of conjecturally finite dimension, with formal properties similar to those of $\Hhat^n(M, p)$, and these are the groups needed in the second author's conjecture on $\zeta$ and $L$-values \cite{Scholbach:specialL}.
%
%
%
%
%
%
%
%
%
\xrema

\rema
In \refth{comparisongroups}, we show that $\Hhat^n(X)$ agrees with $\Khat^T_{-n}(X)$ for $n \leq -1$ and is a subgroup of the latter for $n=0$. The group $\Hhat^1(X) = \coker (K_0(X) \r \oplus \HD^{2p}(X, p))$ is related to the Hodge conjecture, which for any smooth projective $X / \CC$ asserts the surjectivity of $K_0(X)_\Q \r \HD^{2p}(X, \Q(p))$ (Deligne cohomology with rational coefficients). For $n \geq 2$, $\Hhat^n(X) = \oplus \HD^{2p+n-1}(X, p)$.
\xrema


\exam
We list the  groups $\Hhat^{-n} := \Hhat^{-n}(\SpecOF)$ of a number ring $\OF$. These groups and their relation to the Dedekind $\zeta$-function are well-known, cf. \cite[III.4]{Soule:Arakelov}, \cite[p. 623]{Takeda:Higher}.
%
%
%
%
%
For any $n \in \Z$, \refeq{longexactK} reads
$$\HD^0(X, n+1) \r \Hhat^{-2n-1} \r K_{2n+1} \stackrel {\rho_*}\r
\HD^1(X, n+1)
\r \Hhat^{-2n} \r K_{2n}
\stackrel {\rho_*}\r
\HD^0(X, n).$$
In \refth{comparisonRegulator}, we show that the map $\rho_*$ induced by the $\BGL$-module structure of $\oplus \HD \twi p$ agrees with the Beilinson regulator. We conclude by Borel's work that $\Hhat^{-2n-1}$ is an extension of $(K_{2n+1})_{\text{tor}} (=\mu_F \text{ if }n=0)$ by $\HD^0(X, n+1)$ for $n \geq 0$. Moreover, for $n > 0$, $\Hhat^{-2n}$ is an extension of the finite group $K_{2n}$ by a torus, i.e., a group of the form $\RR^{s_n} / \Z^{s_n}$ for some $s_n$ that can be read off \refeq{longExactDeligneBettideRham}. Finally, $\Hhat^0$ is an extension of the class group of $F$ by a group $\RR^{r_1 + r_2 - 1} / \Z^{r_1 + r_2 -1} \oplus \RR$.

For higher-dimensional varieties, the situation is less well-understood. For example, by Beilinson's, Bloch's, and Deninger's work we know that
$$K_{2n+2}(E)_\RR^{(n+2)} \r \HD^2(E, n+ 2)$$
is surjective for $n \geq 0$, where $E$ is a regular proper model of certain elliptic curves over a number field (for example a curve over $\Q$ with complex multiplication in case $n=0$). We refer to \cite[Section 8]{Nekovar:Beilinson's} for references and further examples.
\xexam

\subsection{Functoriality} \mylabel{sect_functoriality}

Let $f: X \r Y$ be a map of $S$-schemes. The structural maps of $X/S$ and $Y/S$ are denoted $x$ and $y$, respectively. We establish the expected functoriality properties of Arakelov motivic cohomology. To define pullback and pushforward, we apply $\Hom_{\DMBei}(-, \HBhatS)$ to appropriate maps, using \refeq{HhatModule}.  

\lemm \mylabel{lemm_pullback}
There is a functorial pullback
$$f^*: \Hhat^n(Y, p) \r \Hhat^n(X, p), \ f^*: \Hhat^n(Y) \r \Hhat^n(X).$$
More generally, for any map $\phi: M \r M'$ in $\SH(S)$ there is a functorial pullback
$$\phi^*: \Hhat^n(M', p) \r \Hhat^n(M, p), \ \ \phi^*: \Hhat^n(M') \r \Hhat^n(M).$$
This pullback is compatible with the long exact sequence \refeq{longExactArakelov} and, for compact objects $M$ and $M'$, with the Arakelov-Chern class \refeq{HBhatVsBGLhat}.
\xlemm
\pf
The second statement is clear from the definition. The first claim follows by applying the natural transformation
$$x_! x^! = y_! f_! f^! y^! \stackrel{\text{\refeq{adjshriek}}}\lr y_! y^!$$
to $\BGL_S$ or $\HBS$, respectively. The last statement is also clear since \refeq{adjshriek} is functorial, in particular it respects the isomorphism $\widehat{\ch}: \BGLhat_{\Q, S} \cong \oplus \HBeihat{S} \twi p$.
\xpf

In the remainder of this section, we assume that $f$ and $y$ (hence also $x$) is a regular projective map (\refde{regularprojectivemap}). Recall that $\dim f = \dim X - \dim Y$.

\defilemm \mylabel{defilemm_pushforward}
We define the pushforward
$$f_*: \Hhat^n(X, p) \r \Hhat^{n-2 \dim(f)}(Y, p-\dim(f))$$
on Arakelov motivic cohomology to be the map induced by the composition
\eqnarr
\M_S(Y) = y_! y^! \HBei{S} & \stackrel{(\tr^\Beilinson_y)^{-1}}\lr & y_! y^* \HBei{S}\twi{\dim(y)} \\
& \stackrel{\text{\refeq{adjstar}}}\lr & y_! f_* f^* y^* \HBei{S}\twi{\dim(y)} \\
& = & x_! x^* \HBei{S}\twi{\dim(y)} \\
& \stackrel{\tr^\Beilinson_x}\lr & x_! x^! \HBei{S}\twi{\dim(y) - \dim(x)} \\
& = & \M_S(X)\twi{-\dim(f)}.
\xeqnarr
Similarly,
$$f_*: \Hhat^n(X) \r \Hhat^{n}(Y)$$
is defined using the trace maps on $\BGL$ instead of the ones for $\HB$ \refeq{traceBGL}, \refeq{traceBeilinson}.

This definition is functorial (with respect to the composition of regular projective maps).
\xdefilemm

\pf
Let $g: Y \r Z$ be a second map of $S$-schemes such that both $g$ (hence $h := g \circ f$) and the structural map $z: Z \r S$ is regular projective. The functoriality of the pushforward is implied by the fact that the following two compositions agree (we do not write $\HBei{-}\twi{-}$ or $\BGL_{-}$ for space reasons):
$$z_! z^! \stackrel{\tr_z^{-1}}\r z_! z^* \r z_! h_* h^* z^* = x_! x^* \stackrel{\tr_x}\r x_! x^!$$
$$z_! z^! \stackrel{\tr_z^{-1}}\r z_! z^* \r z_! g_* g^* z^* = y_! y^* \stackrel{\tr_y}\r y_! y^! \stackrel{\tr_y^{-1}}\r y_! y^* \r y_! f_* f^* y^* = x_! x^* \stackrel{\tr_x}\r x_! x^!$$
This agreement is an instance of the identity $ad_h = y_* ad_f y^* \circ ad_g$.
\xpf

\subsection{Purity and an arithmetic Riemann-Roch theorem} \mylabel{sect_purity}

In this subsection, we establish a purity isomorphism and a Riemann-Roch theorem for Arakelov motivic cohomology. We cannot prove it in the expected full generality of regular projective maps, but need some smoothness assumption.

Given any closed immersion $i: Z \r \SpecZ$, we let $j: U \r \SpecZ$ be its open complement. The generic point is denoted $\eta: \SpecQ \r \SpecZ$. We also write $i$, $j$, $\eta$ for the pullback of these maps to any scheme, e.g., $i: X_Z := X \x_{\SpecZ} Z \r X$. Recall that $B$ is an arithmetic ring whose generic fiber $B_\eta$ is a field (\refde{arithmring}).

Let $f: X \r S$ be a map of regular $B$-schemes. For clarity, we write $\Deligne(p)_{X_\eta}$ for the complex of presheaves on $\Sm/X_\eta$ that was denoted $\Deligne(p)$ above and $\HDel{X_\eta}$ for the resulting spectrum. Moreover, we write $\HDel{X} := \eta_* \HDel{X_\eta} \in \SH(X)$. The complex $\Deligne(p)_{X_\eta}$ is the restriction of the complex $\Deligne(p)_{B_\eta}$. Therefore, there is a natural map $f^* \Deligne(p)_S \r \Deligne(p)_X,$ which in turn gives rise to a map of spectra
$$\alpha_\Deligne^f: f^* \HDel{S} \r \HDel{X}.$$
This map is an isomorphism if $f$ is smooth, since $f^* : \PSh(\Sm/S) \r \PSh(\Sm / X)$ is just the restriction in this case. Is $\alpha_\Deligne^f$ an isomorphism for a closed immersion $f$ between flat regular $B$-schemes? The corresponding fact for $\BGL$, i.e., the isomorphism $f^* \BGL_S = \BGL_X$ ultimately relies on the fact that algebraic $K$-theory of smooth schemes over $S$ is represented in $\SH(S)$ by the infinite Grassmannian, which is a smooth scheme over $S$. Therefore, it would be interesting to have a geometric description of the spectrum representing Deligne cohomology, as opposed to the merely cohomological representation given by the complexes $\Deligne(p)$.

\lemm \mylabel{lemm_preparation}
\begin{enumerate}[(i)]
\item
\mylabel{item_lemmi}
Given another map $g: Y \r X$ of regular $B$-schemes, there is a natural isomorphism of functors $\alpha_\Deligne^g \circ g^* \alpha_\Deligne^f = \alpha_\Deligne^{f \circ g}$.
\item
\mylabel{item_lemmii} The following are equivalent
\begin{itemize}
\item
$\alpha_\Deligne^f$ is an isomorphism in $\SH(X)$.
\item
For any $i: Z \r \SpecZ$, the object $i^! f^* \HDel{S}$ is zero in $\SH(X \x_{\Z} Z)$.
\item
For any sufficiently small $j: U \r \SpecZ$, the adjunction morphism $f^* \HDel{S} \r j_* j^* f^* \HDel{S}$ is an isomorphism in $\SH(X)$.
\end{itemize}
\item \mylabel{item_lemmiii}
The conditions in \refit{lemmii} are satisfied if $f$ fits into a diagram
$$\xymatrix{
X \ar[d]_f \ar[r]^{x}
&
B' \ar[r]
&
B \\
S \ar[ur]_s}$$
where $B'$ is regular and of finite type over $B$, $x$ and $s$ are smooth. In particular, this applies when $f$ is smooth or when both $X$ and $S$ are smooth over $B$.
\end{enumerate}
\xlemm
\pf
\refit{lemmi} is easy to verify using the definition of the pullback functor. \refit{lemmiii} is a consequence of the above remark and \refit{lemmi} using the chain of natural isomorphisms $f^* \HDel{S} = f^* s^* \HDel{B'} = x^* \HDel{B'} = \HDel{X}$.
For \refit{lemmii}, consider the map of distinguished localization triangles
$$\xymatrix{
i_* i^! f^* \HDel{S} \ar[r] \ar[d]
&
f^* \HDel{S} \ar[r] \ar[d]^{\alpha_\Deligne^f}
&
j_* j^* f^* \HDel{S} \ar[d]^{j_* j^* \alpha_\Deligne^f = j_* \alpha_\Deligne^{f_U}} \\
0 = i_* i^! \HDel{X} \ar[r] & \HDel{X} \ar[r] & j_* j^* \HDel{X}.
}
$$
The map $\alpha_\Deligne^{f_U}$ is an isomorphism as soon as $j$ is small enough so that $X_U$ and $S_U$ are smooth over $B_U$. Such a $j$ exists by the regularity of $X$ and $S$. This shows the equivalence of the three statements in \refit{lemmii}.
\xpf

Below, we write $\Beilinson := \oplus_{p \in \Z} \HB \twi p$ and $\widehat \Beilinson_X := \hofib (\Beilinson_X \r \Beilinson_X \wedge \HDel{X})$. We define
$$f^? \BGLhat_S := \hofib (f^! \BGL_S \stackrel {\id \wedge 1} \lr f^! \BGL_S \wedge f^* \HDel{S})$$
and similarly for $f^? \widehat \Beilinson_S$. (The notation is not meant to suggest a functor $f^?$, it is just a shorthand.)  The Chern class $\ch: \BGL_S \r \Beilinson_S$ induces a map $f^? \widehat \ch: f^? \BGLhatS \r f^? \widehat \Beilinson_S$.

\theo \mylabel{theo_absolutepurity}
Let $f: X \r S$ be a regular projective map (\refde{regularprojectivemap}) such that $\alpha^f_\Deligne$ is an isomorphism. (In particular (\refle{preparation} \refit{lemmiii}) this applies when $B$ is a field or when $X$ and $S$ are smooth over $B$ or when $f$ is smooth.) Then there is a commutative diagram in $\SH(X)_\Q$ as follows. Its top row horizontal maps are $\BGL_X$-linear (i.e., induced by maps in $\DMBGL(X)$) and the bottom horizontal maps are $\Beilinson_X$-linear. All maps in this diagram are isomorphisms (in $\SH(X)_\Q$).  
\eqn
\mylabel{eqn_fancy}
\xymatrix{
\BGLhat_X \ar[d]^{\widehat {\ch}_X} &
f^* \BGLhat_S \ar[l]_{\widehat \alpha} \ar[d]^{f^* \widehat{\ch}_S} \ar[rr]^{\widehat{\tr_\BGL}}  &
&
f^? \BGLhat_S \ar[d]^{f^? \widehat {\ch}_S} \ar[r]^\beta &
f^! \BGLhat_S \ar[d]^{f^! \widehat {\ch}_S} \\
\widehat \Beilinson_X &
f^* \widehat \Beilinson_S \ar[l]^{\widehat \alpha} \ar[r]_{\widehat{\Td(T_f)}} &
f^* \widehat \Beilinson_S \ar[r]_{\widehat{\tr_\Beilinson}} &
f^? \widehat \Beilinson_S \ar[r]_\beta &
f^! \widehat \Beilinson_S .
}
\xeqn
\xtheo

\pf
To define the maps $\widehat \alpha$ in \refeq{fancy}, we don't make use of the assumption on $\alpha^f_\Deligne$. Pick fibrant-cofibrant representatives of $\BGL$ and $\HB$, and $\HD$. Thus, in the following diagram of spectra, $f^*$ and $\wedge$ are the usual, non-derived functors for spectra:
$$\xymatrix{
f^* \BGL_S \ar@{=}[r] \ar[d]^{f^*(\id \wedge 1_\Deligne)} &
f^* \BGL_S \ar[r]^{\alpha_\BGL^f} \ar[d]^{\id \wedge f^* 1_\Deligne}&
\BGL_X \ar[d]^{\id \wedge 1_\Deligne} \\
f^* (\BGL_S \wedge \HDel{S}) &
f^* \BGL_S \wedge f^* \HDel{S} \ar[l] \ar[r]^{\alpha_\BGL^f \wedge \alpha_\Deligne^f} &
\BGL_X \wedge \HDel{X}.
}
$$
As $f^*$ is a monoidal functor (on the level of spectra), the canonical lower left hand map is an isomorphism of spectra and the left square commutes. The right square commutes because of $\alpha_\Deligne^f (f^* 1_\Deligne) = 1_\Deligne$. This diagram induces a map between the homotopy fibers of the two vertical maps, which are $f^* \BGLhat_S$ and $\BGLhat_X$, respectively. This is the map $\widehat \alpha$ above. The one for $\widehat \Beilinson$ is constructed the same way by replacing $\BGL$ by $\Beilinson$ throughout. Using $f^* \ch_S = \ch_X$, this shows the commutativity of the left hand square in \refeq{fancy}. By definition of $\BGL$, $\alpha_\BGL^f: f^* \BGL_S \r \BGL_X$ is a weak equivalence.  Thus, both maps $\widehat \alpha$ are isomorphisms in $\SH(X)$ when $\alpha_\Deligne^f$ is so. They are clearly $\BGL_X$- and $\Beilinson_X$-linear, respectively.

The horizontal maps in the middle quadrangle are defined as in \refth{hats}\refit{hats1}: for example, the map $\tr_\BGL : f^* \BGL \r f^! \BGL$ gives rise to $\widehat {\tr_\BGL}: f^* \BGLhat_S \r f^? \BGLhat_S$. It is $\BGL_X$-linear since $\tr_\BGL$ is so. Similarly, we define $\widehat {\Td (T_f)}$ (viewing $\Td(T_f)$ as a ($\Beilinson_X$-linear) map $f^* \Beilinson_S \r f^* \Beilinson_S$) and $\widehat {\tr_\Beilinson}$. Picking representatives of all maps, the quadrangle will in general not commute in the category of spectra, but does so up to homotopy, by construction and by the Riemann-Roch theorem \ref{theo_RiemannRoch}. This settles the middle rectangle.

By the regularity of $X$ and $S$, we can choose $j : U \subset \Spec \Z$ such that $X_U$ and $S_U$ are smooth over $B_U$. We will also write $j$ for $X_U \r X$ etc. 

By assumption, $\alpha^f_\Deligne$ is an isomorphism. Hence, the adjunction map $f^! \BGL \wedge f^* \HD \r j_* j^* (f^! \BGL \wedge f^* \HD)$ is an isomorphism in $\SH$. In fact, both terms are isomorphic in $\SH$ to $\oplus_p \HD \twi p$, as one checks for example using the purity isomorphism $f^! \BGL_S \cong f^* \BGL_S = \BGL_X$. Thus, $f^? \BGLhat$ is canonically isomorphic to the homotopy fiber of $f^! \BGL \r j_* j^* f^! \BGL \r j_* j^* (f^! \BGL \wedge f^* \HD) = j_* (j^* f^! \BGL \wedge j^* f^* \HD)$. Here, the last equality is a canonical isomorphism on the level of spectra, since $j^*$ is just the restriction. By definition, $j^* f^! = j^! f^!$. We may therefore replace $f$ by $f_U$. Now, $f_U^! M$ is \emph{functorially} isomorphic (in $\SH$) to $f_U^* M \twi{n}$, $n := \dim f_U$, by construction of the relative purity isomorphism by Ayoub \cite[Section 1.6]{Ayoub:Six1}. Indeed, $a$ is a closed immersion,  and $p$ and every map in the diagram with codomain $B_U$ is smooth:
$$\xymatrix{
X_U \ar[r]_a \ar@/^/[rr]^{f_U} \ar[drr] &
\P[n]_{S_U} \ar[dr] \ar[r]_p &
S_U \ar[d] \\
&
&
B_U.
}$$

Hence to construct $\beta$, it is enough to construct a commutative diagram of spectra
$$\xymatrix{
f_U^* \BGL_{S_U} \twi n \ar[d]^{\id \wedge 1_\Deligne} \ar@{=}[r] &
f_U^* \BGL_{S_U} \twi n \ar[d]^{f^*(\id \wedge 1_\Deligne)} \\
f_U^* \BGL_{S_U} \twi n \wedge f_U^* \HDel{S_U} \ar[r]^\gamma &
f_U^* (\BGL_{S_U} \wedge \HDel{S_U}) \twi n.
}
$$
The map $\gamma$ is the natural map of spectra $f_U^* x \wedge f_U^* y \r f_U^* (x \wedge y)$, which clearly makes the diagram commute in the category of spectra. We have constructed a map (in $\SH$) $\beta: f^? \BGLhat \r f^! \BGLhat$ in a way that is functorial with respect to (a lift to the category of spectra of) the map $\ch : \BGL \r \Beilinson$. Therefore, the analogous construction for $\widehat \Beilinson$ produces the desired commutative square of isomorphisms (in $\SH$). Again, the top row map $\beta$ is $\BGL$-linear and the bottom one is $\Beilinson_X$-linear.

Finally, the vertical maps in \refeq{fancy} are isomorphisms using the Arakelov Chern character \refeq{ArakelovChernCharacter}. 
\xpf

We can now conclude a higher arithmetic Riemann-Roch theorem. It controls the failure of $\widehat \ch$ to commute with the pushforward.
\theo \mylabel{theo_HARR}
Let $f: X \r S$ be a regular projective map (\refde{regularprojectivemap}) of schemes of finite type over an arithmetic ring $B$ (\refde{arithmring}). Moreover, we assume that $f$ is such that
$$\alpha^f_\Deligne : f^* \HDel{S} \r \HDel{X}$$
is an isomorphism. This condition is satisfied, for example, when $f$ is smooth or when $X$ and $S$ are smooth over $B$ (\refle{preparation}). Then, the following holds:

\begin{enumerate}[(i)]
\item \mylabel{item_purity} (Purity)
The absolute purity isomorphisms for $\BGL$ and $\HB$ \refeq{absolutepurity} induce isomorphisms (of $\BGL_X$- and $\HB{}_X$-modules, respectively): 
$$\BGLhat_X \cong f^* \BGLhat_S \cong f^! \BGLhat_S, \ \ \HBhat_X \cong f^* \HBhat_S \cong f^! \HBhat_S \twi {-\dim f}.$$
In particular, Arakelov motivic cohomology is independent of the base scheme in the sense that there are isomorphisms
$$\Hhat^n(X / S) \cong \Hhat^n (X / X), \ \ \Hhat^n(X / S, p) \cong \Hhat^n (X / X, p).$$
\item \mylabel{item_HARR} (Higher arithmetic Riemann-Roch theorem)
There is a commutative diagram
$$\xymatrix{
\Hhat^n(X / X) \ar[r]^{f_*} \ar[d]^{\widehat \ch_X} & \Hhat^n(S/S) \ar[d]^{\widehat \ch_S} \\
\oplus_{p \in \Z} \Hhat^{n+2p}(X, p) \ar[r]^{f_* \circ \widehat{\Td}(T_f)} & \oplus_{p \in \Z} \Hhat^{n+2p}(S, p).
}$$
Here, the top line map $f_*$ is given by
\begin{eqnarray*}
\Hhat^n(X / X) & := & \Hom_{\SH(X)} (S^{-n}, \BGLhat_X)
\\
& \stackrel{\text{\refeq{fancy}}}\r &  \Hom_{\SH(X)}(S^{-n}, f^! \BGLhat_S) 
\\
& \stackrel{\text{\refeq{adjstar}}} \r &  \Hom_{\SH(S)}(S^{-n}, \BGLhat_S) = \Hhat^n (S / S).
\end{eqnarray*}
Using the identifications $\Hhat^n(X/X) \cong \Hhat^n(X/S)$, this map agrees with the one defined in \ref{defilemm_pushforward}. The bottom map $f_*$ is given similarly replacing $\BGL$ with $\Beilinson$.
\end{enumerate}
\xtheo

\pf
The isomorphisms for $\BGLhat_?$ in \refit{purity} are a restatement of \refth{absolutepurity}. The ones for $\HBhat_?$ also follow from that by dropping the isomorphism $\widehat {\Td(T_f)}$ in the bottom row of \refeq{fancy} and noting that $\tr_\Beilinson$, hence $\widehat{\tr_\Beilinson}$ shifts the degree by $\dim f$. The isomorphisms in the second statement are given by the following identifications of morphisms in $\DMBGL(-)$, using \refeq{HhatModule}:
\begin{eqnarray*}
\Hom(\BGL_X, \BGLhat_X) & \stackrel{\text{\ref{theo_absolutepurity}}} \lr & \Hom (\BGL_X, f^! \BGLhat_S) \\
& \stackrel{(\tr_\BGL)^{-1}} \lr & \Hom(f^! \BGL_S, f^! \BGLhatS) \\
&= &\Hom (f_! f^! \BGL_S, \BGLhatS).
\end{eqnarray*}
and likewise for $\HB$.

\refit{HARR} is an immediate corollary of \refth{absolutepurity}, given that the two isomorphisms (in $\SH(X)_\Q$) $\widehat {\Td(T_f)} \circ \widehat \alpha^{-1}$ and $\widehat \alpha^{-1} \circ \widehat {\Td(T_f)}$, where $\Td(T_f)$ is seen as an endomorphism of $f^* \Beilinson_S$ and of $\Beilinson_X$, respectively, agree. This agreement follows from the definition of $\widehat \alpha$. The agreement of the two definitions of $f_*$ is clear from the definition.
\xpf

This also elucidates the behavior of \refeq{longexactK} with respect to pushforward: in the situation of the theorem, the pushforward $f_*: \Hhat^n(X) \r \Hhat^n(S)$ sits between the usual $K$-theoretic pushforward and the pushforward on Deligne cohomology (which is given by integration of differential forms along the fibers in case $f(\CC)$ is smooth, and by pushing down currents in general), multiplied with the Todd class (in Deligne cohomology) of the relative tangent bundle.

\subsection{Further properties} \mylabel{sect_further}
\theo \mylabel{theo_immediate}
\begin{enumerate}[(i)]
\item \mylabel{item_immediatei}
Arakelov motivic cohomology satisfies h-descent (thus, a fortiori, Nisnevich, \'etale, cdh, qfh and proper descent). For example, there is an exact sequence
$$\dots \r \Hhat^n(X, p) \r \Hhat^n(U \sqcup V, p) \r \Hhat^n(W, p) \r \Hhat^{n+1}(X, p) \r \dots$$
where
$$\xymatrix{W \ar[r] \ar[d] & V \ar[d]^p \\
U \ar[r]^f & X}$$
is a cartesian square of smooth schemes over $S$ that is either a distinguished square for the cdh-topology ($f$ a closed immersion, $p$ is proper such that $p$ is an isomorphism over $X \backslash U$), or a distinguished square for the Nisnevich topology ($f$ an open immersion, $p$ \'etale inducing an isomorphism $(p^{-1}(X \backslash U)_\red \r (X \backslash U)_\red$), or such that $U \sqcup V \r X$ is an open cover.
\item \mylabel{item_immediateii}
Arakelov motivic cohomology is homotopy invariant and satisfies a projective bundle formula:
$$\Hhat^n(X \x \A, p) \cong \Hhat^n(X, p),$$
$$\Hhat^n(\mathbf P(E), p) \cong \oplus_{i=0}^d \Hhat^{n-2i}(X, p-i).$$
Here $X / S$ is arbitrary (of finite type), $E \r X$ is a vector bundle of rank $d+1$, $\mathbf P(E)$ is its projectivization.
\item \mylabel{item_immediateiii}
Any distinguished triangle of motives induces long exact sequences of Arakelov motivic cohomology. For example, let $X/S$ be an l.c.i.\ scheme (\refex{lci}). Let $i: Z \subset X$ be a closed immersion of regular schemes of constant codimension $c$ with open complement $j: U \subset X$. Then there is an exact sequence
$$\Hhat^{n-2c}(Z, p-c) \stackrel{i_*} \r \Hhat^n(X, p) \stackrel{j^*}\r \Hhat^n(U, p) \r \Hhat^{n+1-2c}(Z, p-c).$$
\item
The cdh-descent and the properties \refit{immediateii}, \refit{immediateiii} hold \emph{mutatis mutandis} for $\Hhat^*(-)$.
\end{enumerate}
\xtheo
\pf
The h-descent is a general property of modules over $\HBS$ \cite[Thm 16.1.3]{CisinskiDeglise:Triangulated}. The $\A$-invariance and the bundle formula are immediate from $\M(X) \cong \M(X\x \A)$ and $\M(\mathbf P(E)) \cong \oplus_{i=0}^d \M(X)\twi{i}$. For the last statement, we use the localization exact triangle \cite[2.3.5]{CisinskiDeglise:Triangulated} for $U \stackrel j \r X \stackrel i \leftarrow Z$:
$$f_! j_! j^! f^! \HBS \r f_! f^! \HBS \r f_! i_* i^* f^! \HBS.$$
The purity isomorphism $f^* \HBS \twi{\dim f} = f^! \HBS$ (\refex{lci}) for the structural map $f: X \r S$ and the absolute purity isomorphism \refeq{absolutepurity} for $i$ imply that the rightmost term is isomorphic to $f_! i_! i^! f^! \HBS \twi{-\dim i} = \M_S(Z) \twi{-c}$. Mapping this triangle into $\HBhatS(p)[n]$ gives the desired long exact sequence.

The arguments for $\BGLhatS$ are the same. The only difference is that descent for topologies exceeding the cdh-topology requires rational coefficients.
\xpf

\bibliography{bib}
\end{document}